\documentclass[12pt,a4paper]{article}
\pdfoutput=1
\usepackage{booktabs}
\usepackage{amsmath,amssymb,amsfonts}
\usepackage{algorithm}
\usepackage{array}
\usepackage{algorithmicx}
\usepackage{algpseudocode}
\usepackage{graphicx}
\usepackage{diagbox}
\usepackage{textcomp}
\usepackage{mathrsfs}
\usepackage{hhline}
\usepackage{setspace}
\usepackage{color}
\usepackage{soul}
\usepackage{eucal}
\usepackage{float}
\usepackage{subfigure}
\usepackage{indentfirst}
\usepackage{enumerate}
\usepackage[capitalize]{cleveref}
\usepackage{epstopdf}
\usepackage{sidecap}
\usepackage{caption}
\usepackage{multirow}
\usepackage{xcolor}
\usepackage{colortbl,booktabs}

\newtheorem{theorem}{Theorem}[section]

\newtheorem{definition}{Definition}[section]
\newtheorem{example}{Example}[section]
\newtheorem{Remark}{Remark}[section]

\def\BibTeX{{\rm B\kern-.05em{\sc i\kern-.025em b}\kern-.08em
		T\kern-.1667em\lower.7ex\hbox{E}\kern-.125emX}}

\allowdisplaybreaks

\begin{document}

\title{\bf A novel class of fractional adams method for solving uncertain fractional differential equation}

\author{Chenlei Tian$^{1}$, Jing Cao$^1$,~Yifu Song$^1$, Ting Jin$^1$\thanks{Corresponding author.}   \\
{\small\em $^1$School of Science, Nanjing Forestry University}\\
{\small\em Nanjing 210037, China}\\
{\small\em Email: targaryentcl@njfu.edu.cn (C. Tian), jing@njfu.edu.cn (J. Cao)}\\
{\small\em 	songyifu@njfu.edu.cn(Y. Song), tingjin@njfu.edu.cn (T. Jin)}}
\date{}
\maketitle

\begin{abstract}
Uncertain fractional differential equation (UFDE) is a kind of differential equation about uncertain process. As an significant  mathematical tool to describe the evolution process of dynamic system, UFDE is better than the ordinary differential equation with integer derivatives because of its  hereditability and memorability characteristics. However, in most instances, the precise analytical solutions of UFDE is difficult to obtain due to the complex form of the UFDE itself. Up to now, there is not plenty of researches about the numerical method of UFDE, as for the existing numerical algorithms, their accuracy is also not high. In this research, derive from the interval weighting  method, a class of fractional adams method is innovatively proposed to solve UFDE. Meanwhile, such fractional adams method extends the traditional predictor-corrector method to higher order cases. The stability and truncation error limit of the improved algorithm are analyzed and deduced. As the application, several numerical simulations (including $\alpha$-path, extreme value and the first hitting time of the UFDE) are provided to manifest the higher accuracy and efficiency of the proposed numerical method.
\end{abstract}

\section{Introduction}

The most common mathematical tools used to describe uncertain phenomena are frequency-based probability theory and degree of belief based on uncertainty theory. In real life,  the precondition of using probability theory to deal with random phenomena is that there is sufficient sample data so that the probability distribution can be inferred to approximate the actual situation. Then the probability theory can solve many problems and really have good effects. Owing to practical or technical reasons, sufficient sample data cannot be obtained in most cases, so the evaluation must be based on the trust of experts in the field. It is worth mentioning that the Nobel Prize winner in economics Kahneman and  Tversky \cite{kahneman2013prospect} pointed out that humans generally overestimate the likelihood of some events happening, which can make the results deviate from reality and even make decision makers make wrong decisions. To workaround the belief degree of experts mathematically, Liu \cite{liu2007uncertainty,baoding2007uncertainty} established the uncertainty theory in 2007. In fact, the uncertainty theory just describes the case where the confidence level is usually much greater than the cumulative frequency for small or unavailable sample size. For more details about uncertainty theory, see reference \cite{liu2012there}.

After putting forward the uncertainty theory, Liu \cite{liu2007uncertainty}  defined the uncertain measure, uncertain space, uncertain variable together with its uncertain distribution, expectation and variance respectively. For more details, please refer to \cite{liu2009some}. Liu \cite{liu2008fuzzy} defined the uncertain process for the sake of describing the evolution of uncertain phenomenon over time in 2008. In order to study the uncertain calculus of uncertain process, Liu \cite{liu2009some} proposed Liu process, which is a Lipschitz continuous uncertain process with normal uncertain variables as stationary independent increments. Moreover, he advanced the uncertain differential equations (UDEs) to describe the evolution of uncertainty over time. After that, researchers began to pay attention to the numerical method of solving UDE because UDE has many applications in uncertain optimal control \cite{zhu2010uncertain} and uncertain finance \cite{chen2013uncertain,yang2020parameter,liu2013toward}. In 2013, Yao and Chen \cite{yao2013numerical} first putted forward the famous Yao-Chen formula.The relationship between ordinary differential equation and UDE has thus been established. Based on the Yao-Chen formula, a series of numerical methods have emerged: Among them, Yang and Ralescu \cite{yang2015adams} proposed the Adams method to solve the UDE, Milne method was designed by Gao\cite{gao2016milne} in 2016 for solving the UDEs.  Wang and Ning \cite{wang2015adams} provided Adams-Simpson method, while Zhang and Gao \cite{zhang2017hamming} proposed Hamming method to solve the UDEs.

\begin{table}[H]
	\renewcommand\arraystretch{1.5}
	\centering
	\caption{Comparison of numerical approaches for uncertain (fractional) differential equations}
	\begin{tabular}{m{1.0cm}<{\centering}m{3.7cm}<{\centering}m{1.9cm}<{\centering}m{3.5cm}<{\centering}m{4cm}<{\centering}}
		\toprule
		Work & Author & Order &    Algorithm &  Solution \\
		\midrule
		\cite{yao2013numerical} & Kai Yao \&Xiaowei Chen &  $p=1$  & Euler's method &  Analytical solution \&Numerical solution \\
		\midrule
		\cite{yang2015adams} & Xiangfeng Yang \&Dan A Ralescu & $p=1$ &  Adams method &  Numerical solution \\
		\midrule
		\cite{gao2016milne} &  Rong Gao & $p=1$ &  Miline methodd &  Numerical solution \\
		\midrule
		\cite{wang2015adams} & Xiao Wang \&Yufu Ning & $p=1$ & Adams-simpson method &  Numerical solution \\
		\midrule
		\cite{zhang2017hamming} & Yi Zhang etal. & $p=1$ & Haming method & Analytical solution \\
		\midrule
		\cite{Wu2022sest} & GuoCheng Wu et al. & $p\in(0,1]$ & Adams method & Numerical solution \\
		\midrule
		\cite{Luo_Wu_Huang_2022} & Luo Cheng et al. & $p\in(0,1]$ & Truncation method & Analytical solution \&Numerical solution \\
		\midrule
		\cite{lu2019numerical} & Ziqiang Lu \&Yuanguo  Zhu & $p\in(0,1]$ & Truncation method & Numerical solution \\
		\bottomrule
	\end{tabular}
\end{table}

In some complex dynamic systems, ordinary differential equations usually cannot adequately describe the complex operating mechanisms and hereditability of the system. In this case, fractional differential equation (FDE) can be well used for describing the memorability and historical features. The fractional calculus method has been widely used in the financial field in recent years \cite{ma2020application,liping2021new}, optimal control \cite{kilbas2006theory} and image encryption \cite{wu2019new,abdeljawad2020discrete}. For the latest researches on FDEs, please refer to \cite{wu2022unified,fu2021fractional}. In 2013, Zhu \cite{zhu2015uncertain} first combined the fractional theory with the uncertainty theory, defined two types of UFDEs, the Caputo one and the Riemann-Liouville one. Subsequently, based on linear growth condition and Lipschitz condition, Zhu \cite{zhu2015existence} proposed the existence and uniqueness theorem of the UFDE's solution. Using the $\alpha$-path in Yao-Chen formula, the relationship between the FDEs and the UFDEs is established, the solution of UFDEs can be represented by a cluster of solutions of FDEs, that is, $\alpha$-path is a numerical methods for solving UFDE. Lu \cite{lu2019numerical} proposed a numerical approach to solve UFDE involving Caputo derivatives by using $\alpha$-path, and also provided a formula to calculate the expected value of the monotone function of the solution about UFDE. In many other fields, UFDE also has its own applications. Jin \cite{jin2019extreme} studied the extreme value of a kind of Caputo type UFDE's solution and subsequently applied it to the American option pricing model. In addition, Jin \cite{jin2020first} also used the predictor-corrector method to give the uncertain distributions of their first hitting time of a kind of nonlinear Caputo type UFDE, and applied it to a novel uncertain risk index model. Wu et al. \cite{Wu2022sest} studied the parameter estimation of UFDE based on fractional Adams method. Considering the generalization of UFDE, Luo et al. \cite{Luo_Wu_Huang_2022} studied the uniqueness and existence of the solution of the generalized fractional uncertain differential equation (GUFDE), and gives the extreme values and solutions of GUFDEs. It is worth noting that in these studies, the accuracy of the method used to calculate the uncertain distribution is not too high, so it is necessary to study a kind of numerical algorithm with higher accuracy and faster arithmetic speed.

As the main motivation, our work aims to study the improvement of fractional order Adams numerical algorithm. Kai Diethelm \cite{diethelm2002predictor} first proposed the fractional Adams method. Li \cite{li2011numerical} further proposed the fractional Adams method based on Simpson method, which approximated the fractional integral and the fractional derivative with the use of higher order piecewise interpolating polynomials. This approach was designed to solve the analytical expression of the uncertain integral. Thus when there are too many nodes, the complexity of the solution will be very large. Based on the motivation of improving the computational accuracy of Adams method, our study proposes a novel fractional Adams method which can be extended to any node.

The frame of article is as follows: In Section 2, some concepts and properties of the UFDE, coupled with the product integration method and Adams method, are reviewed. In Section 3, the Adams method is extended to order $n$. Relevant numerical experiments can be found in Section 4, which calculate the extreme value, the inverse distribution and the first hitting time (FHT) of UFDE, respectively. Finally, a brief conclusion is given in Section 5.

\section{Preliminary}

In this section, some definitions and theorems about uncertainty theory are introduced. The brief introductions of the Adams method and UFDE are given.

\subsection{Uncertainty theory and UFDE}

\begin{definition}(Liu\cite{baoding2007uncertainty})
Assume that $\Gamma$ is a nonempty set, $\mathcal{L}$ is a $\sigma$-algebra over $\Gamma$. Call each element $\Lambda\in \mathcal{L}$ an event. Call a set function $\mathcal{M}$ defined on the $\sigma$-algebra an uncertain measure, if the following four axioms are satisfied:

\quad a. $\mathcal{M}\{\Gamma\}=1$.

\quad b. $\mathcal{M}\{\Lambda\}+\mathcal{M}\left\{\Lambda^{c}\right\}=1$ for any event $\Lambda$.

\quad c. $\mathcal{M}\left\{\cup_{i=1}^{\infty} \Lambda_{i}\right\} \leq \sum\limits_{i=1}^{\infty} \mathcal{M}\left\{\Lambda_{i}\right\}$ for events $\Lambda_{1}, \Lambda_{2}, \ldots$

\quad d. Let $\left(\Gamma_{k}, \mathcal{L}_{k}, \mathcal{M}_{k}\right)~(k\in\mathbb{N}^+)$ be uncertain spaces. The product uncertain measure $\mathcal{M}$ is an uncertain measure satisfied
$$\mathcal{M}\left\{\prod_{k=1}^{\infty} \Lambda_{k}\right\}=\bigwedge\limits_{k=1}^{\infty} \mathcal{M}_{k}\left\{\Lambda_{k}\right\}$$.
\end{definition}

\begin{definition}(Liu\cite{liu2009some})
Call an uncertain process $C_t$ a canonical Liu process, if

\quad a. $C_0=0$ and almost all sample paths are Lipschitz continuous.

\quad b. $C_t$ has stationary and independent increments.

\quad c. every increment $C_{s+t}-C_s$ is a normal uncertain variable with expected value $0$ and variance $t^2$.

The uncertain distribution of $C_t$ is
\begin{equation*}\label{Phi_distribution}
\Phi_t(x)=\left(1+\exp\left(-\frac{\pi x}{\sqrt{3}t}\right)\right)^{-1}
\end{equation*}
and inverse uncertain distribution is
\begin{equation*}\label{Phi_inv}
\Phi_t^{-1}(\alpha)=\frac{\sqrt{3}t}{\pi}\ln\frac{\alpha}{1-\alpha}.
\end{equation*}
\end{definition}


\begin{definition}(Yao-Chen\cite{yao2013numerical})
An UDE
\begin{equation*}\label{UDE}
\mathrm{d} X_t=f\left(t, X_t\right) \mathrm{d} t+g\left(t, X_t\right) \mathrm{d} C_t
\end{equation*}
is said to have an $\alpha$-path $X_t^\alpha$ if it solves the  differential equation
\begin{equation*}
\mathrm{d} X_t^\alpha=f\left(t, X_t^\alpha\right) \mathrm{d} t+\left|g\left(t, X_t^\alpha\right)\right| \Phi^{-1}(\alpha) \mathrm{d} t,~0<\alpha<1
\end{equation*}
where
\begin{equation*}
\Phi^{-1}(\alpha)=\frac{\sqrt{3}}{\pi}\ln\frac{\alpha}{1-\alpha}.
\end{equation*}
is the inverse standard normal uncertain distribution.
\end{definition}


\begin{theorem}(Yao-Chen\cite{yao2013numerical})
Let $X_t$ be the solution, $X_t^\alpha$ be the $\alpha$-path of the UDE\eqref{UDE}, respectively. Thus

$$\mathcal{M}\left\{X_t \leq X_t^\alpha, \forall t\right\}=\alpha,$$
$$\mathcal{M}\left\{X_t>X_t^\alpha, \forall t\right\}=1- \alpha.$$
\end{theorem}

\begin{definition}(\cite{Podlubny1999FractionalDE})
For any $\nu>0, t>a$, the Riemann-Liouville type fractional integral is defined as
\begin{equation*}
I^\nu_{a}f(t)=\frac{1}{\Gamma(\nu)}\int^t_a(t-s)^{\nu-1}f(s,X_s)\mathrm ds,
\end{equation*}
where $\Gamma(\alpha)$ is Gamma Function.
\end{definition}


\begin{definition}(\cite{Podlubny1999FractionalDE})
Assume that $\nu>0, t>a$ and $n-1<\nu\leq n~(n\in\mathbb N^+$), Caputo type fractional differential is defined as
\begin{equation*}
^C_aD_t^\nu f(t)=I^{n-\nu}_{a+}f^{(n)}(t)=\frac{1}{\Gamma(n-\nu)}\int^t_a(t-s)^{n-\nu-1}f^{(n)}(s,X_s)\mathrm ds.
\end{equation*}
\end{definition}


\begin{definition}(Li\cite{li2011numerical})
Suppose $C_t$ is a canonical Liu process, $\nu>0$, $f:[0,+\infty)\times\mathbb R^n\to\mathbb R^n$ and $g:[0,+\infty)\times\mathbb R^n\to\mathbb R^{n+l}$ are two given functions. An UFDE
\begin{equation*}
_a^CD_t^\nu X_t=f\left(t, X_t\right)+g\left(t, X_t\right) \frac{\mathrm{d} C_t}{\mathrm{d} t.}
\end{equation*}

The corresponding $\alpha$-path $X_t^\alpha$ of $X_t$ is a function of $t$ which solves the following FDE
\begin{equation*}
_a^C D_t^{\nu} X_t^\alpha=f\left(t, X_t^\alpha\right)+\left|g(t,X_t^\alpha)\right|\Phi^{-1}(\alpha).
\end{equation*}
\end{definition}


\begin{theorem}(Zhu\cite{zhu2015existence})
An UFDE has a unique solution $X_t$ in $[0,+\infty)$, if the coefficients $f(t,X_t)$ and
$g(t,X_t)$ satisfy the Lipschitz condition
\begin{equation*}
|f(t, x)-f(t, y)|+|g(t, x)-g(t, y)| \leq L|x-y|, \quad \forall x, y \in \mathbb R^n, \quad t \in[0,+\infty)
\end{equation*}
and the linear growth condition
\begin{equation*}
|f(t, x)|+|g(t, x)| \leq L(1+|x|), \quad \forall x \in \mathbb R^n, \quad t \in[0,+\infty).
\end{equation*}
\end{theorem}


\begin{theorem}(Lu\cite{lu2019european})
A Caputo type UFDE with initial value conditions
\begin{equation}\label{UFDE_init}
\left\{\begin{array}{l}
\displaystyle{ }_a^C D_t^{\nu} X_t=f\left(t, X_t\right)+g(t,X_t)\frac{\mathrm dC_t}{\mathrm dt}\\
\left.X_t^{(i)}\right|_{t=a}=x_i, i=0,1, \ldots, n-1
\end{array}\right.
\end{equation}
has an integral form

\begin{equation*}\label{UFDE_intform}
X_t=\sum_{i=0}^{n-1}\frac{(t-a)^ix_i}{\Gamma(i+1)}+\frac{1}{\Gamma(\nu)} \int_a^t(t-s)^{\nu-1}f\left(s, X_s\right) \mathrm{d} s+\frac{1}{\Gamma(\nu)} \int_a^t(t-s)^{\nu-1}g\left(s, X_s\right) \mathrm{d} C_s,
\end{equation*}
and if its $\alpha$-path $X_t^\alpha$ satisfies

\begin{equation}\label{UFDE_init_inv}
\left\{\begin{array}{l}
\displaystyle{ }_a^C D_t^{\nu} X_t^\alpha=f\left(t, X_t^\alpha\right)+g(t,X_t^\alpha)\frac{\mathrm dC_t}{\mathrm dt},~   \alpha\in(0,1) \\
\left.[X_t^\alpha]^{(i)}\right|_{t=a}=x^\alpha_i, i=0,1, \ldots, n-1,
\end{array}\right.
\end{equation}
then
\begin{equation}\label{UFDE_inv_intform}
\begin{aligned}
X_t^\alpha =&\sum_{i=0}^{n-1} \frac{(t-a)^i x_i^\alpha}{\Gamma(i+1)}+\frac{1}{\Gamma(v)} \int_a^t(t-s)^{v-1} f\left(s, X_s^\alpha\right) \mathrm{d} s \\
&+\frac{1}{\Gamma(v)} \int_a^t(t-s)^{v-1}\left|g\left(s, X_s^\alpha\right)\right| \Phi^{-1}(\alpha) \mathrm{d} s.
\end{aligned}
\end{equation}
\end{theorem}

\section{Improved $n$-order fractional Adams method}
 In this section, the numerical method of a class of integral equations, namely the product integration method, and the  fractional Adams method, will be introduced. On this basis, this paper will improve the Fractional Adams Method to apply it to any number of nodes. This paper proves the operability of this method and gives the truncation error.
\subsection{Improved fractional Adams method}

Product Integration Method was originally proposed to solve the integral equation. The numerical solution of fractional differential equation can be converted into the numerical method of the second kind of Volterra integral equation with a circular singular kernel. The problem of The fractional Adams method proposed by Kai Diethelm \cite{li2011numerical} is the application of two-point Lagrange interpolation in Product Integration Method.

We can take
\begin{equation}\label{s_k}
\hat{I}=\sum\limits_{i=0}^{n}{\sum\limits_{j=0}^{m}{{{w}_{ij}}}}f ({{t}_{ij}}),
\end{equation}
as an approximate estimate of integral
\begin{equation*}
I=\int^b_ap(s)f(s)\mathrm ds
\end{equation*}
where ${{w}_{ij}}=\int_{t_{i}}^{t_{i+1}}p(s)L_{ij}(s)\mathrm{d}t$ and $L_{ij}$ is base function of Lagrange interpolation.

Let us consider the integral equation with $0<\nu\leq1$
\begin{equation*}
X(t_{n+1})=x_0+\frac{1}{\Gamma(\nu)}\int^{t_{n+1}}_{t_0}(t_{n+1}-s)^{\nu-1}f(s,X_s)\mathrm ds,
\end{equation*}
using fractional Adams method, we have
\begin{equation*}
\int^{t_{n+1}}_{t_0}(t_{n+1}-s)^{\nu-1}f(s,X_s)\mathrm ds\approx\sum\limits_{j=0}^{n+1}a_{j,n+1}f(t_j,X_{t_j})
\end{equation*}
where $a_{j,n+1}=\int^{t_{n+1}}_{t_0}(t_{n+1}-s)^{\nu-1}\mathcal L_{j,n+1}(s)\mathrm ds$ and
\begin{equation*}
\mathcal{L}_{j, n+1}(s)= \begin{cases}\left(s-t_{j-1}\right) /\left(t_j-t_{j-1}\right) & , \text { if } t_{j-1}<s<t_j \\ \left(t_{j+1}-s\right) /\left(t_{j+1}-t_j\right) & , \text { if } t_j<s<t_{j+1} \\ 0 & , \text { otherwise }\end{cases}.
\end{equation*}

In the premise of stability, this paper tries to explore a higher order fractional Adams method, which means increasing the number of nodes of Lagrange interpolation. Li \cite{li2011numerical} proposed the fractional Adams method based on Simpson method in 2011. However, solving the coefficient before the node involves solving the equation
\begin{equation*}
\int^{t_{i+1}}_{t_i}(t_{n+1}-s)^{\nu-1}L_{ij}(s)\mathrm ds
\end{equation*}\
or the corresponding analytical expression of indefinite integral, with the number of nodes increases, the difficulty of solving is also increasing. Therefore, we need to improve the fractional Adams method to obtain a more accurate solution.

Back to the fractional Adams method, if we want to use more points to improve the calculation accuracy, we need to give a general formula to calculate the term $\displaystyle\int_{t_i}^{t_{i+1}}(t_{n+1}-s)^{\nu-1}s^k\mathrm ds$. Only when this term has an exact value or analytical formula in different cases of $k$ can it be applied to fractional integration. Now we give the calculation method of the exact value of this item.

\begin{theorem}\label{theorem3.1}
Let $M(\nu-1,k)=\displaystyle\int_{t_i}^{t_{i+1}}(t_{n+1}-s)^{\nu-1}s^k\mathrm ds$, then
\begin{equation}
M(\nu-1,k)=-\frac{{{s}^{k}}}{k+\nu }{{(t_{n+1}-s)}^{\nu }}\bigg|_{t_i}^{t_{i+1}}+\frac{kt_{n+1}}{k+\nu }M(\nu-1,k-1).	
\end{equation}
\end{theorem}
{\bf Proof:} For $\displaystyle\int_{t_i}^{t_{i+1}}(t_{n+1}-s)^{\nu-1}s^k\mathrm ds$, we have
\begin{equation*}
\int_{t_i}^{t_{i+1}}{{{(t_{n+1}-s)}^{\nu -1}}{{s}^{k}}\mathrm{d}s}=-\frac{{{s}^{k}}}{\nu }{{(t_{n+1}-s)}^{\nu }}\bigg|_{t_i}^{t_{i+1}}+\frac{k}{\nu }\int_{t_i}^{t_{i+1}}{{{(t_{n+1}-s)}^{\nu }}{{s}^{k-1}}\mathrm{d}s},
\end{equation*}
\begin{equation*}
\int_{t_i}^{t_{i+1}}{{{(t_{n+1}-s)}^{\nu }}{{s}^{k-1}}\text{d}s}={t}\int_{t_i}^{t_{i+1}}{{{(t_{n+1}-s)}^{\nu -1}}{{s}^{k-1}}\mathrm{d}s}-\int_{t_i}^{t_{i+1}}{{{(t_{n+1}-s)}^{\nu -1}}{{s}^{k}}\mathrm{d}s}.
\end{equation*}
That is,
\begin{equation}\label{M_1_eq}
M(\nu-1,k)=-\frac{{{s}^{k}}}{\nu }{{(t_{n+1}-s)}^{\nu }}\bigg|_{t_i}^{t_{i+1}}+\frac{k}{\nu }M(\nu ,k-1),
\end{equation}
\begin{equation}\label{M_2_eq}
M(\nu ,k-1)={t}M(\nu -1,k-1)-M(\nu -1,k).
\end{equation}
Simultaneous Eq. (\ref{M_1_eq}) and Eq. (\ref{M_2_eq}), we have
\begin{equation}\label{M_eq}
M(\nu -1,k)=-\frac{{{s}^{k}}}{\nu }{{(t_{n+1}-s)}^{\nu }}\bigg|_{t_i}^{t_{i+1}}+\frac{k}{\nu }\left[ tM(\nu -1,k-1)-M(\nu -1,k) \right].
\end{equation}
Thus, the theorem is proved. 

\begin{Remark}
According to Theorem \ref{theorem3.1}, when the form of the integral $\displaystyle\int_{t_i}^{t_{i+1}}({{t}_{n+1}}-s)^{\nu-1}\mathrm{d}s$ is known, we can easily compute the integral $\displaystyle\int_{t_i}^{t_{i+1}}(t_{n+1}-s)^{\nu-1}s^k\mathrm ds$ for any $k\geq1$. The original complex integral calculation is transformed into solving the relationship between different $k \geq 1$, which greatly reduce the computation complexity.
\end{Remark}

In this way, we successfully find the relation about $\displaystyle\int_{t_i}^{t_{i+1}}(t_{n+1}-s)^{\nu-1}s^k\mathrm ds$ under different $k$. In addition, most node-related data generated during the calculation process can be reused to avoid redundancy.

\begin{theorem}
	The product integral over the interval $[t_i,t_{i+1}]$
	\begin{equation}\label{I_hat}
		\hat{I}=\int_{{{t}_{i}}}^{{t_{i+1}}}(t_{n+1}-s)^{\nu-1}\sum\limits_{j=0}^mf (t_{ij}){L}_{ij}(s)\mathrm{d}s
	\end{equation}
	 can be written as
	\begin{equation}\label{I_hat_lagrange_polynomial}
		\hat I=\sum\limits_{k=0}^{m-1}\left[\sum\limits_{j=0}^mf(t_{ij})a_{jk}\right]\int^{t_{i+1}}_{t_i}(t_{n+1}-s)^{\nu-1}s^k\mathrm ds.
	\end{equation}
	Then, rewrite the Lagrange interpolation for $m$ nodes into a form of polynomial
	\begin{equation*}
		L_{ij}=a_{j0}s^0+a_{j1}s^1+\cdots+a_{jm-1}s^{m-1}.
	\end{equation*}

\end{theorem}
{\bf Proof:}
	Substituting the expression of $L_{ij}$ into $\displaystyle \int_{{{t}_{i}}}^{{{t}_{i+1}}}(t_{n+1}-s)^{\nu-1}{L}_{ij}(s)\mathrm{d}s$, we can obtain
	\begin{equation*}
		\int_{{{t}_{i}}}^{{{t}_{i+1}}}(t_{n+1}-s)^{\nu-1}\sum\limits_{k=0}^{m-1}a_{jk}s^k\mathrm{d}s,
	\end{equation*}
	then substitute it into Eq. (\ref{I_hat}), the conclusion apparently be proofed.

%

In this way, combining the relation of the integral terms with Eq. (\ref{I_hat_lagrange_polynomial}), we can get the expression of the numerical algorithm.

Considering the Lagrange interpolation method for node $\{t_i\}_{0\leq i\leq n}$,

\begin{algorithm}[H]
  \caption{Splitting of Lagrange polynomials.}
  \label{alg:Lagrange}
  \begin{algorithmic}[1]
  \State Set $t=[t_0,t_1,\dots,t_m];$
  \label{code:fram:1}
  \State For $j=0,1,\dots,m$, perform step \ref{code:fram:3} to step \ref{code:fram:7} successively;
  \label{code:fram:2}
  \State $\text{deno}=\prod\limits_{i\not=j}^m(t_i-t_j);$
  \label{code:fram:3}
  \State For $k=[0,1,\dots,m-1]$, perform step \ref{code:fram:5} to step \ref{code:fram:6} successively;
  \label{code:fram:4}
  
  \State Choose $m-k$ values from $[\dots, t_{j-1}, t_{j+1}, \dots]$ to compute the product, there are $C^k_m$ cases, the sum of them is denoted as $a_{jk}$;
  \label{code:fram:5}
  \State $a_{jk}=(-1)^{m-k}a_k;$
  \label{code:fram:6}
  \State $C_j=[a_{j0},a_{j1},\cdots,a_{jm-1}]/\text{deno};$
  \label{code:fram:7}
  \State Output the result $C_j=[a_{j0},a_{i1},\cdots,a_{jm-1}], j=0,1,\dots,m.$
  \end{algorithmic}
\end{algorithm}

 where the formula and polynomial form of the interpolation of the $j$ term are
\begin{equation*}
{{L}_{j}}(s)=\frac{\left( s-{{t}_{0}} \right)\cdots \left( s-{{t}_{j-1}} \right)\left( s-{{t}_{j+1}} \right)\cdots \left( s-{{t}_{m}} \right)}{\left( {{t}_{j}}-{{t}_{0}} \right)\cdots \left( {{t}_{j}}-{{t}_{j-1}} \right)\left( {{t}_{j}}-{{t}_{j+1}} \right)\cdots \left( {{t}_{j}}-{{t}_{m}} \right)}
\end{equation*}
and
$$
{{L}_{j}}(s)=a_0s^{0}+a_1s^1+\cdots +a_{m-1}s^{m-1},
$$
repectively.

\subsection{ Predictor-Corrector method of improved fractional Adams method}
Given the node $t_{k-m},t_{k-m+1},\cdots,t_{k-1}$ and the corresponding value $X_{t_{k-m}},X_{t_{k-m+1}},\cdots,X_{t_{k-1}}$, we estimate 
\begin{equation*}
X_{t_{n}}=X_{t_{n-1}}+\frac{1}{\Gamma(\nu)} \int_{t_{n-1}}^{t_n}(t_n-s)^{\nu-1} f(s, X_s)\mathrm ds
\end{equation*} 
of the interval $[t_{k-1},t_k]$. Let the predictor term be
\begin{equation}
X^p_{t_n}=X_{t_{n-1}}+\frac{1}{\Gamma(\nu)}\sum\limits_{j=n-m}^{n-1}{w}_jf(t_j,X_{t_j}),
\end{equation}
where
\begin{equation}
{{w}_{j}}=\int_{{{t}_{n-1}}}^{t_{n}}(t_n-s)^{v-1}{{L}_{j}}(s)\mathrm{d}s,
\end{equation}
\begin{equation}
L_j(s)=\prod_{\substack{n-m \leq i \leq n-1 \\ i \neq j}} \frac{s-t_i}{s_j-t_i}.
\end{equation}
Then the corrector term is
\begin{equation}
X_{t_n}=X_{t_{n-1}}+\frac{1}{\Gamma(\nu)}\left[\sum\limits_{j=n-m+1}^{n-1}{w}_jf(t_j,X_{t_j})+w_jf(t_n,X_{t_n})\right].
\end{equation}

\begin{theorem}
Considering $f(t)\in C^n[a,t]$, for the $n$-th order fractional Adams method on the subinterval $[t_{n-1},t_n]$, the truncation error limit is:
\begin{equation}
\begin{aligned}
\left|X_t-\hat{X}_t\right| &=\left|I_a^v[f(s)-\hat{f}(s)]\right| \\
& \leq|f(s)-\hat{f}(s)| \frac{1}{\Gamma(v)} \int_{t_{n-1}}^{t_n}(t-s)^{v-1} \mathrm{~d} s \\
&=\frac{\left(t-t_{n-1}\right)^v-\left(t-t_n\right)^v}{\Gamma(v+1)} \frac{M_n}{n !}\left|\omega_n(s)\right|,
\end{aligned}
\end{equation}\label{accuracy_pr}
where $\omega_(s)=(s-t_0)(s-t_1)\cdots(s-t_{n-1}),M_n=\max\limits_{a\leq s\leq t}f^{(n)}(x)$.
\end{theorem}

To solve the FDE as follow
\begin{equation*}
\left\{\begin{array}{l}
\displaystyle{ }_a^C D_t^{\nu} X_t=f\left(t, X_t\right)\\
\left.X_t^{(i)}\right|_{t=a}=x_i, i=0,1, \ldots, n-1.
\end{array}\right.
\end{equation*}
The following algorithm is given.
\begin{algorithm}[H]
  \caption{Predictor-corrector fractional Adams method of order $n$.}
  \label{alg:Adams}
  \begin{algorithmic}[1]
  \State Set interval $t\in[a,s]$, initial value $X_a$, and $n$-order.
  \State Let $t=[0.01:0.01:1], x_{\text{seq}}=X_a, f_{\text{seq}}=[\,\,];$
  \State For $m=[1,\dots,n-1]$, First, the formula Eq. (\ref{M_eq}) and Algorithm \ref{alg:Lagrange} are used to calculate the predictor value $x_{\text{seq}}(m+1)$; then the formula and algorithm are applied to obtain the fractional Adams correction value of order $m+1$ and replace the predicted value.
  \State For $i=[1,2,\dots, \text{length}(t_{\text{seq}})-1]$, use $x_{\text{seq}}(i), t_{\text{seq}}$ compute $n$ order predicted value $x_{\text{seq}}(i+1)$; use $x_{\text{seq}}(i+1), t_{\text{seq}}$ compute $n$ order correction value $x_{\text{seq}}(i+1)$ and replace the predicted value
  \State Output the results $x_{\text{seq}}(end)$
  \end{algorithmic}
\end{algorithm}
In the case of limiting the number of intervals to be divided, it is recommended to divide the intervals first intensively then incompactly, which is because in the case of only one initial value, before using $n$ order Adams algorithm, we need to use $m=1,2,\cdots,n-1$ order Adams algorithm to calculate the first $n$ points successively. Therefore, the error of the initial point will affect the subsequent calculation, so the intervals of the first few points is reduced to avoid its impact on the accuracy in the subsequent weighted calculation. Or you can use a more expensive method to calculate the value of the first few points with higher accuracy, and then use $n$ order Adams algorithm for calculation.

\section{Numerical Simulations}
This section introduces an example of using the modified $n$ order Adams algorithm to calculate the extreme value of UFDE. The differential equation used in this example has the analytical solution. This paper give the error analysis of the inverse distribution $X^{\alpha}_t$ of UFDE's solution for $\alpha \in [0.01,0.99]$ and $t\in[0.01,1]$. The calculation results using the algorithm are given, and the magnitude of error for different $n$ are compared.

Based on two differential equations without analytical solutions, numerical methods are constructed. One of which is applied to the inverse distribution of UFDE's solution, and the other is applied to calculate solve the first hitting time of the UFDE. Both of them have achieved good results.
\subsection{Extreme value of UFDE's solution}

In this subsection, we calculate extreme values of UFDE's solution with a numerical format, which is described in Algorithm \ref{alg:Extreme_value}. According to the description in \cite{jin2019extreme}, let $X_t$ is the unique solution for Eq. (\ref{UFDE_init}) and the $X_t^\alpha$ is $\alpha$-path. $\Psi_t$ is the uncertain distribution of $X_t$. The IUD of the infimum $\inf\limits_{0\leq t\leq s}J(X_t)$ exists,

\begin{equation}\label{eq:ide}
\bar{\Upsilon}_s^{-1}(\alpha)= \begin{cases}\inf\limits_{0 \leq t \leq s} J\left(\Psi_t^{-1}(\alpha)\right) &  \text { if } J(x) \text { increases strictly } \\ \inf\limits_{0 \leq t \leq s} J\left(\Psi_t^{-1}(1-\alpha)\right) &  \text { if } J(x) \text { decreases strictly }\end{cases}.
\end{equation}

\begin{algorithm}[H]
  \caption{Extreme value of UFDE's solution}
  \label{alg:Extreme_value}
  \begin{algorithmic}[1]
	  \State Write UFDE in the corresponding integral form and split $t$ and $\alpha$ according to the interval 0.01.
	  \State Compute the inverse distribution of UFDE's solution about different $t\in[a, b]$ and $\alpha\in(0,1)$ using algorithm \ref{alg:Adams}.
	  \State Calculate the extreme value of the solution for different $\alpha$ according to Eq. (\ref{eq:ide}).
  \end{algorithmic}
\end{algorithm}

\begin{example}
Assume the following linear Caputo type UFDE
\begin{equation}\label{eg1}
_a^C D_t^{\nu}X_t=aX_t+b t^{\upsilon}\frac{\mathrm dC_t}{\mathrm dt},
\end{equation}
with initial value $X_0=0.5$ and $J(x)=x$.
\end{example}

\begin{figure}[H]
	\centering  
	\subfigbottomskip=1pt 
	\subfigcapskip=-5pt 
	\subfigure[$n=2$]{\label{figure1_a}
		\includegraphics[width=0.45\linewidth]{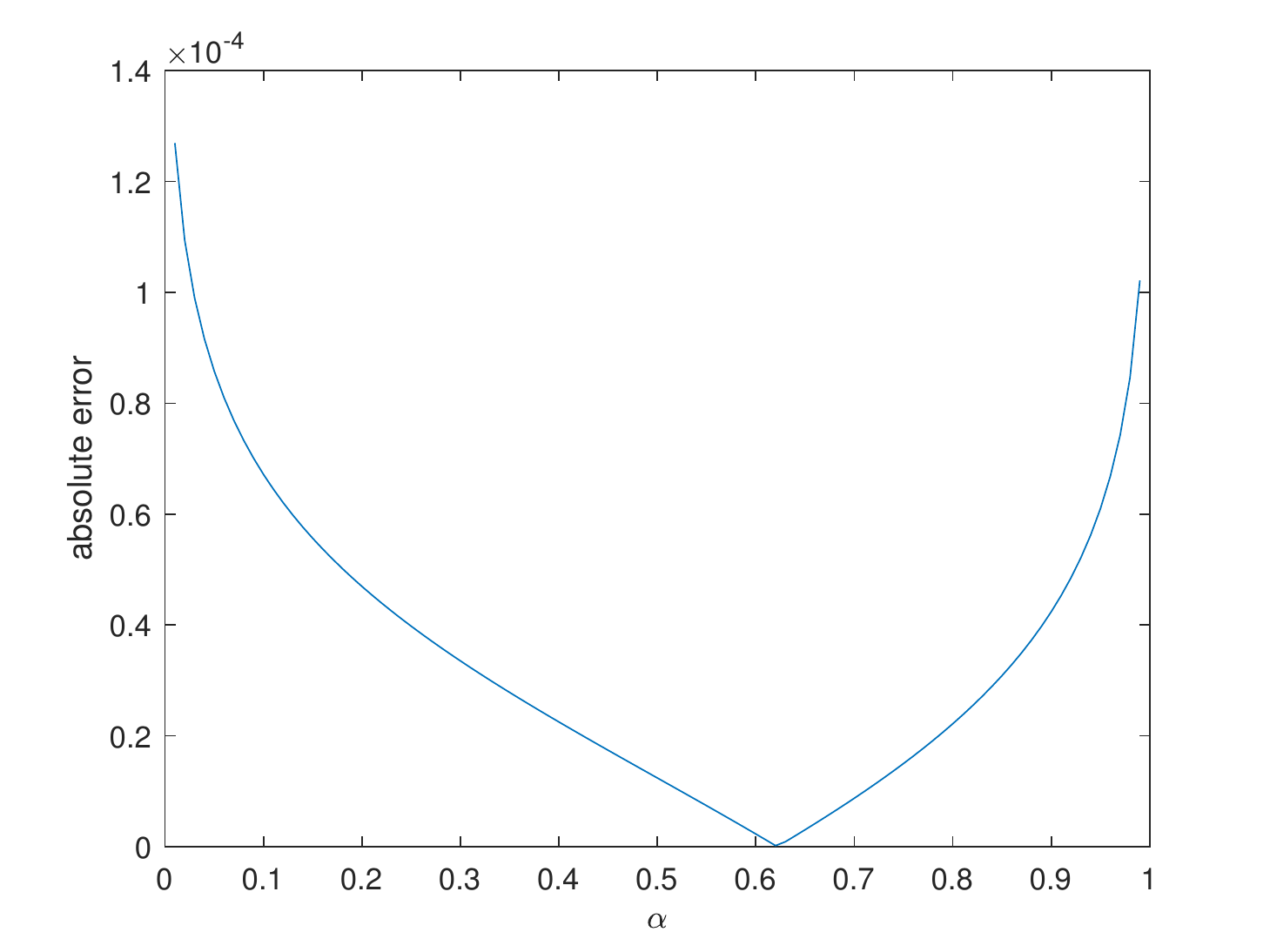}}
	\subfigure[$n=2$]{\label{figure1_d}
		\includegraphics[width=0.45\linewidth]{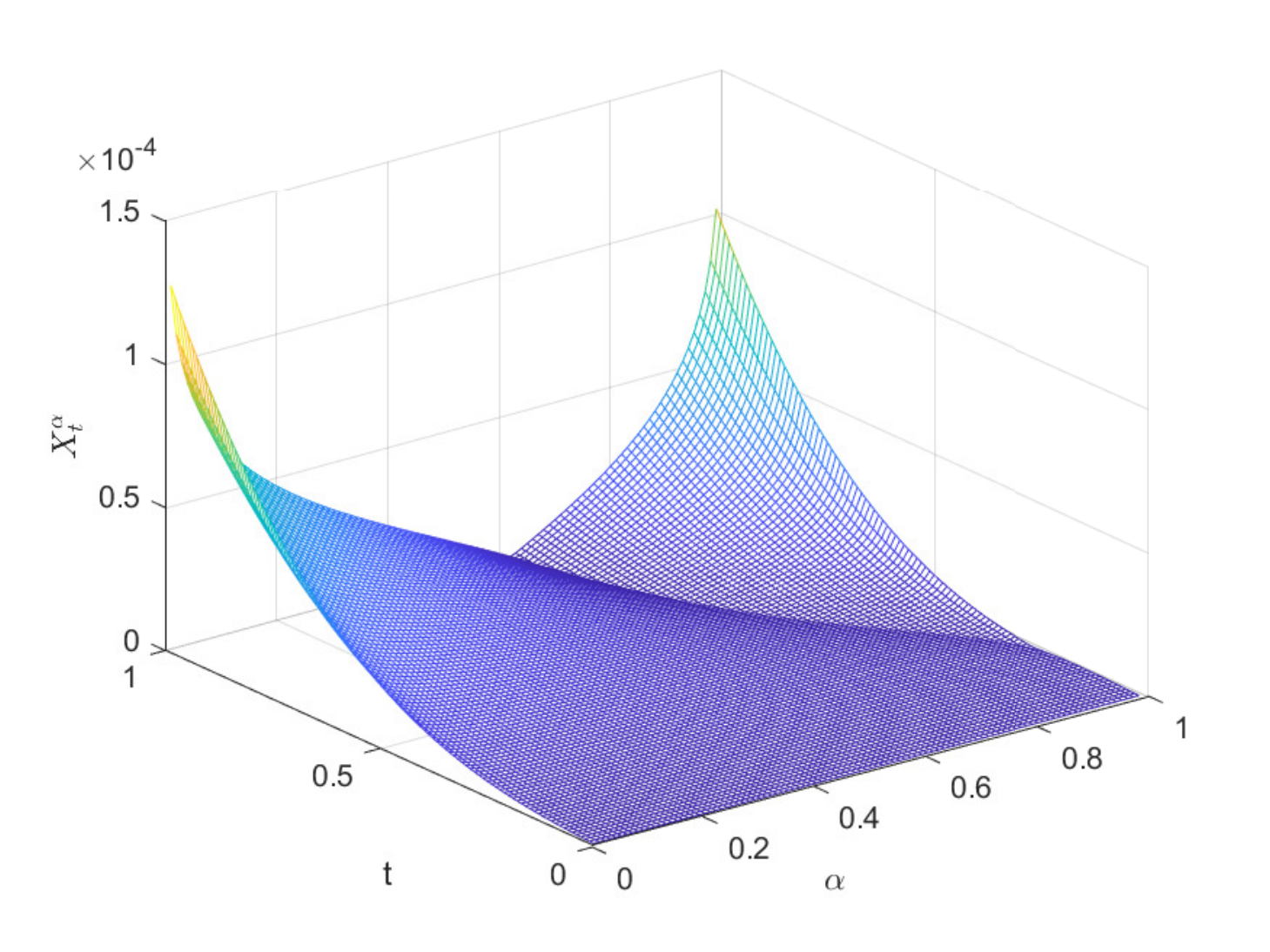}}
\\
	\subfigure[$n=3$]{\label{figure1_b}
		\includegraphics[width=0.45\linewidth]{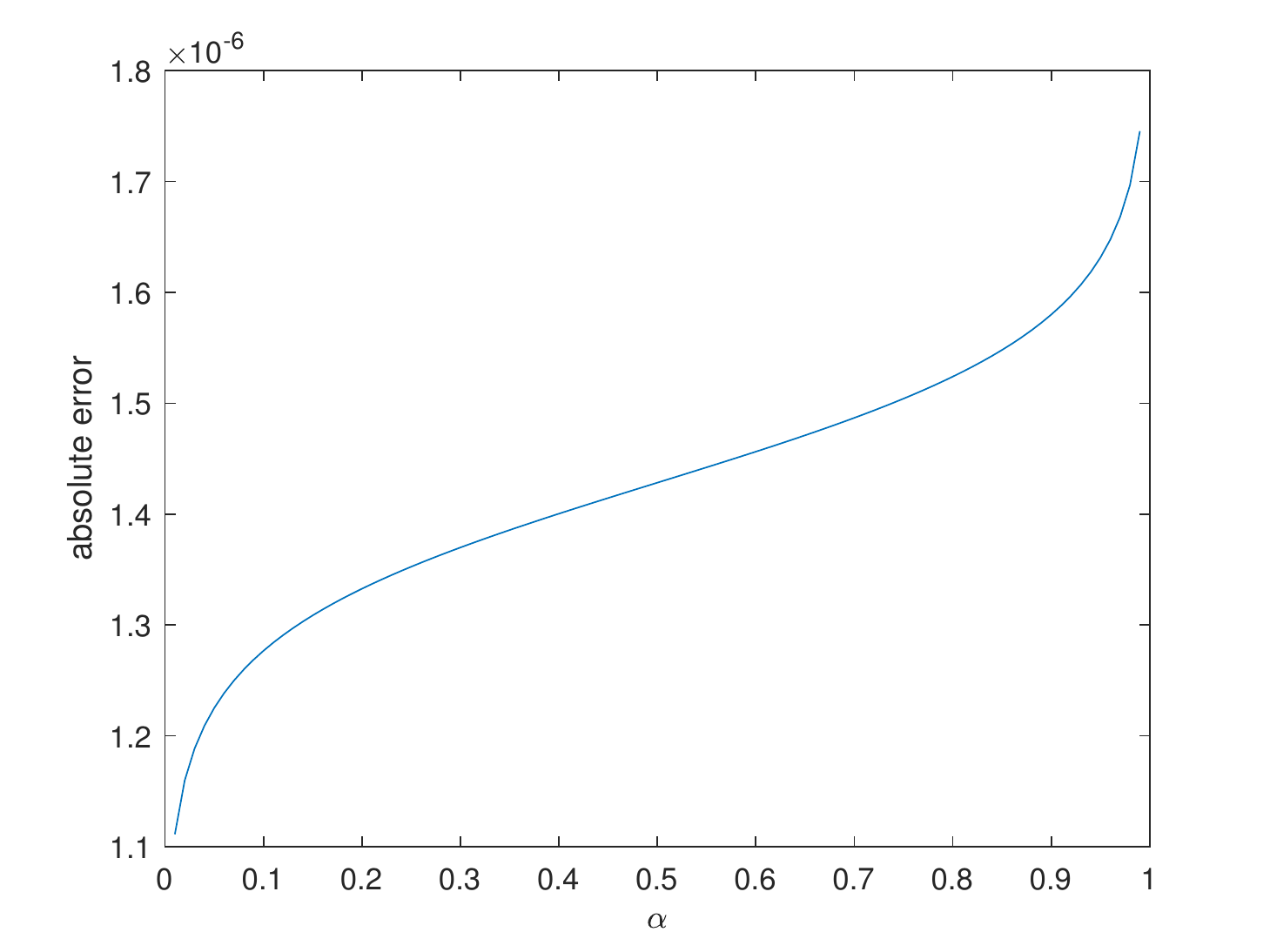}}
	\subfigure[$n=3$]{\label{figure1_e}
		\includegraphics[width=0.45\linewidth]{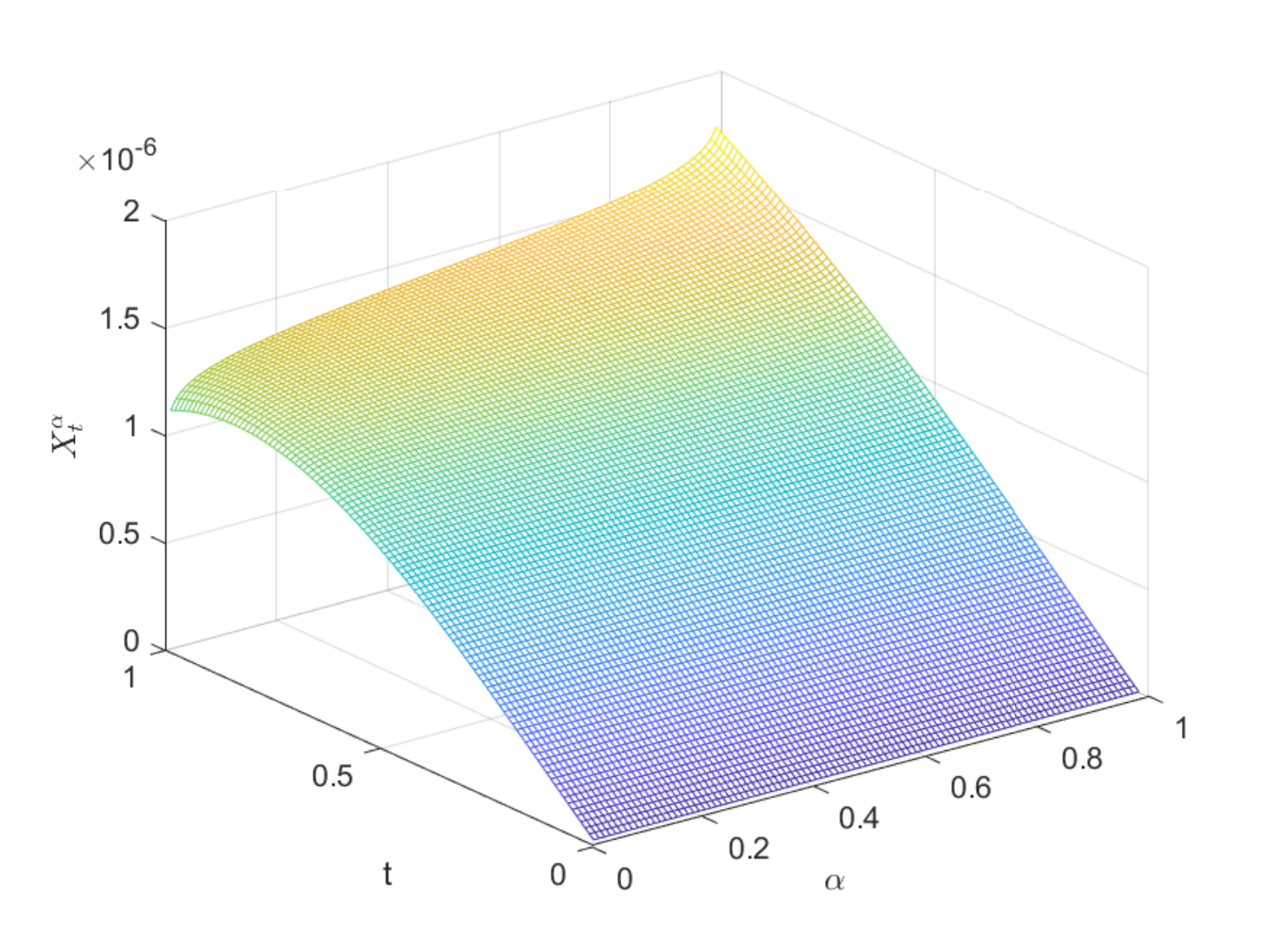}}
\\
	\subfigure[$n=4$]{\label{figure1_c}
		\includegraphics[width=0.45\linewidth]{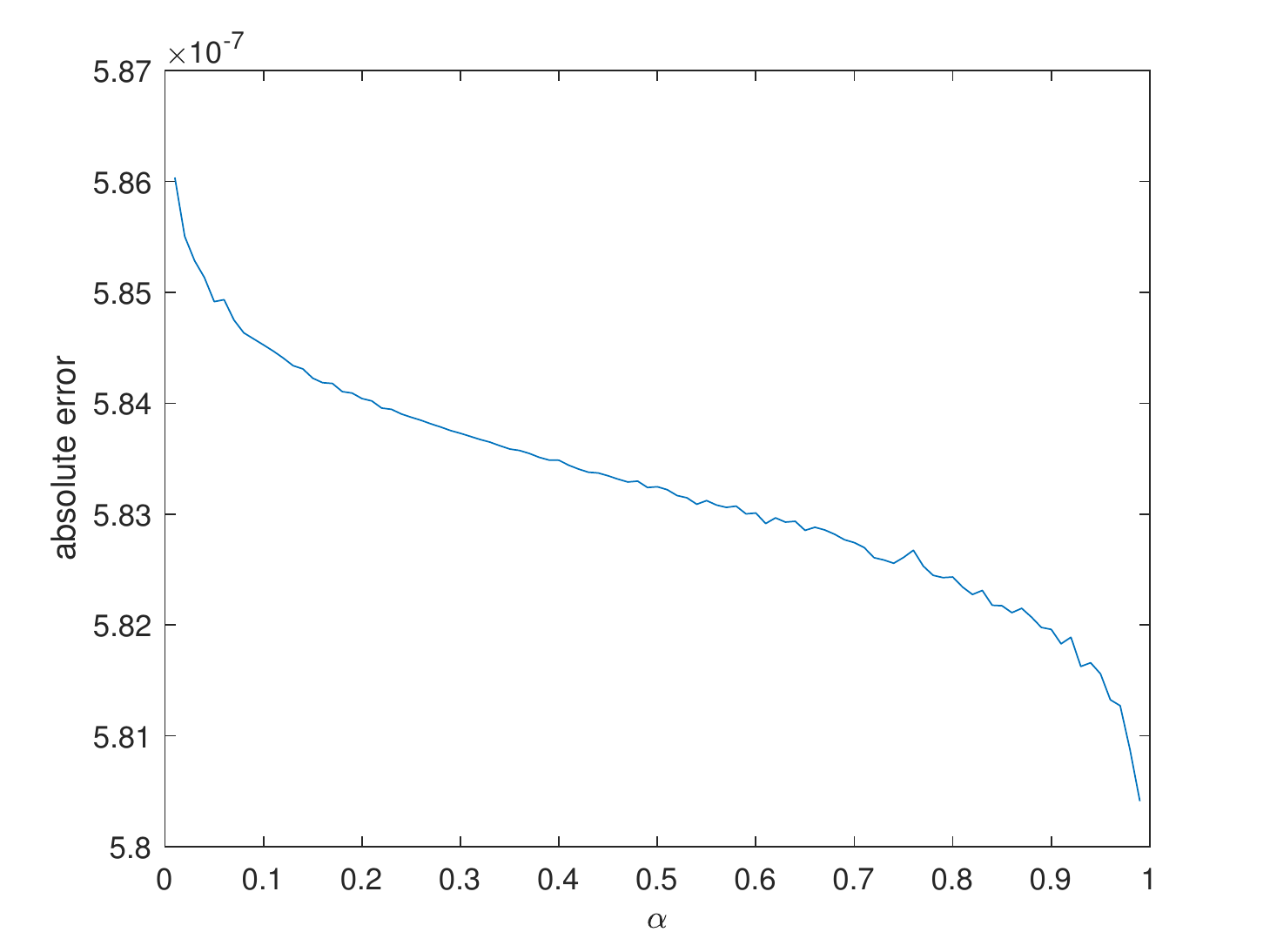}}
	\subfigure[$n=4$]{\label{figure1_f}
		\includegraphics[width=0.45\linewidth]{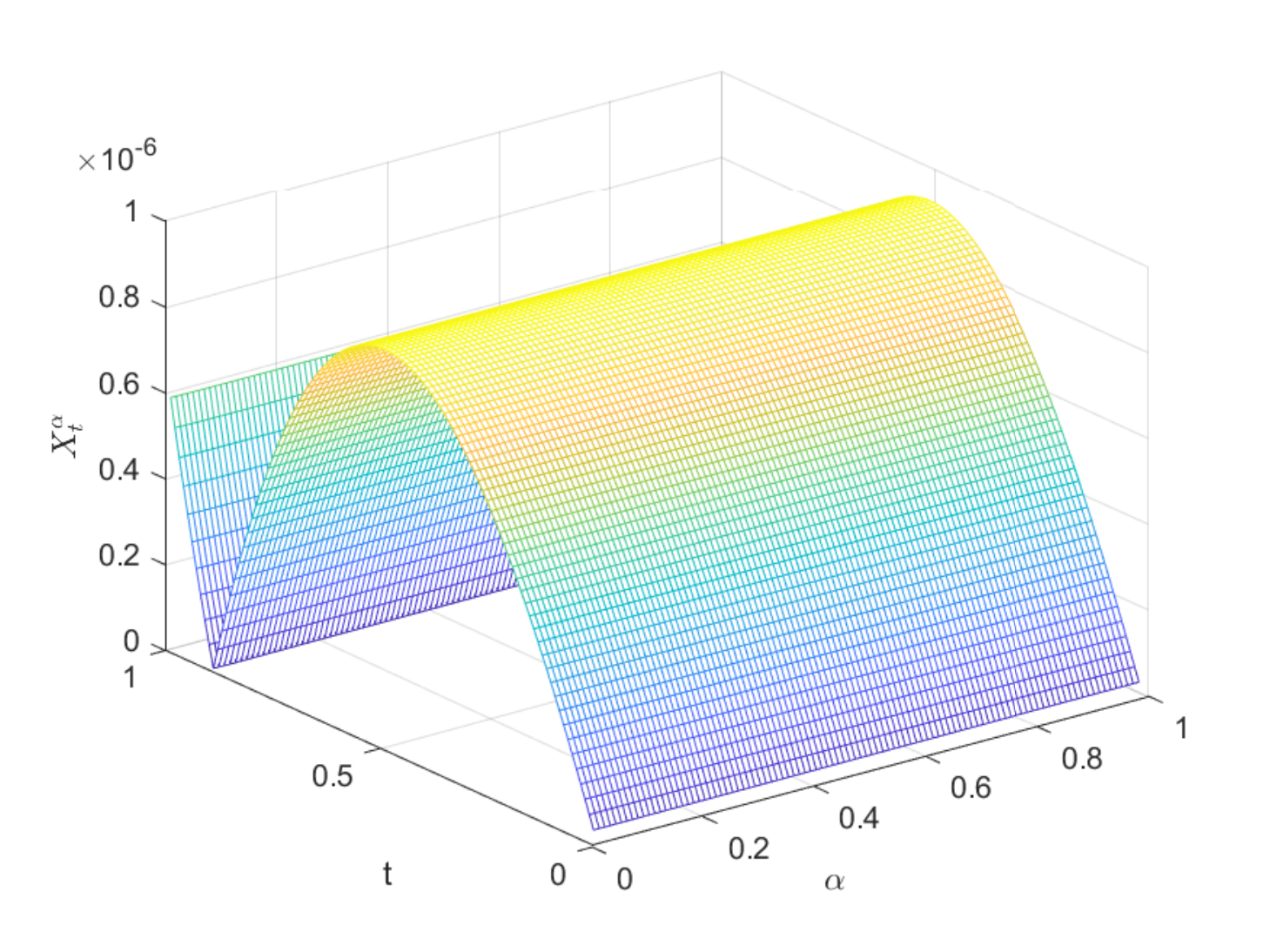}}
	\caption{Plots of absolute error and $\alpha$-paths of (\ref{eg1})}
\end{figure}

We calcuate the absolute errors of different $\alpha$-path with parameter values $a=0.6, b=1, t=1, \upsilon=2$ at $n\in[2,3,4]$, see \cref{figure1_a,figure1_b,figure1_c}. Moreover, the absolute error of uncertain distributions of $X^{\alpha}_t$ with different $t$ and $\alpha$ at $n\in[2,3,4]$ provided by Algorithm \ref{alg:Adams} is shown in \cref{figure1_d,figure1_e,figure1_f}.
The uncertain distribution of extreme value is shown in Fig. \ref{figure2_a} and the absolute error of extreme value is shown in Fig. \ref{figure2_b}. In order to show the absolute error under different $n$ better, we perform a logarithmic on the error results.

\begin{figure}[H]
	\centering  
		\includegraphics[width=0.7\linewidth]{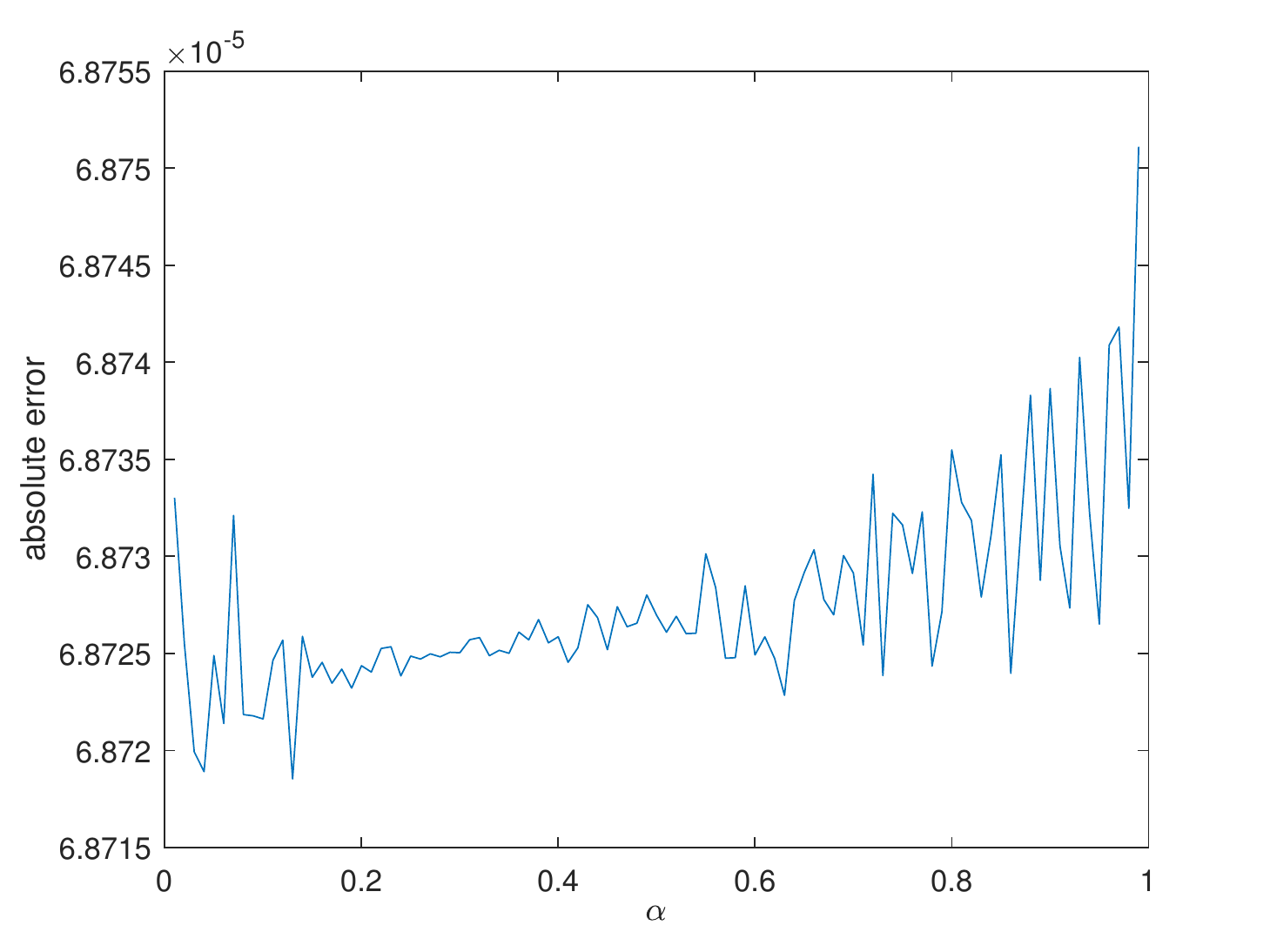}
	\caption{Absolute error about $\alpha$ with $n=5$}\label{figure1sub}
\end{figure}

As can be seen from \cref{figure1_a,figure1_b,figure1_c}, when $n=2,3,4$, the order of error is $10^{-4},\,\,10^{-6}$ and $10^{-7}$. In \cref{figure1_d,figure1_e,figure1_f}, when $n$ increases , the change rate of error about $t$ also increases accordingly with the same value $\alpha$. In fact, as the number of nodes increases, the stability of the algorithm will also deteriorate. In Fig. \ref{figure1_c}, when $n=4$, part of the function image about the error of $\alpha$ appears jagged.

When $n=5$, in Fig. \ref{figure1sub}, this situation worsens. This image of the absolute value of error about $\alpha$ shows that it is very sensitive to the change of the inverse distribution of diffusion term $bt^\nu\Phi^{-1}(\alpha)$. For this example, if $n$ is too large, its effect will not be quite as impressive in terms of sensitivity or accuracy. Considering that the essence of fractional Adams method is the application of Lagrange interpolation, this phenomenon may be caused by over-fitting caused by too many nodes.

\begin{figure}[H]
	\centering  
	\subfigbottomskip=1pt 
	\subfigcapskip=-5pt 
	\subfigure[Extreme value]{\label{figure2_a}
		\includegraphics[width=0.48\linewidth]{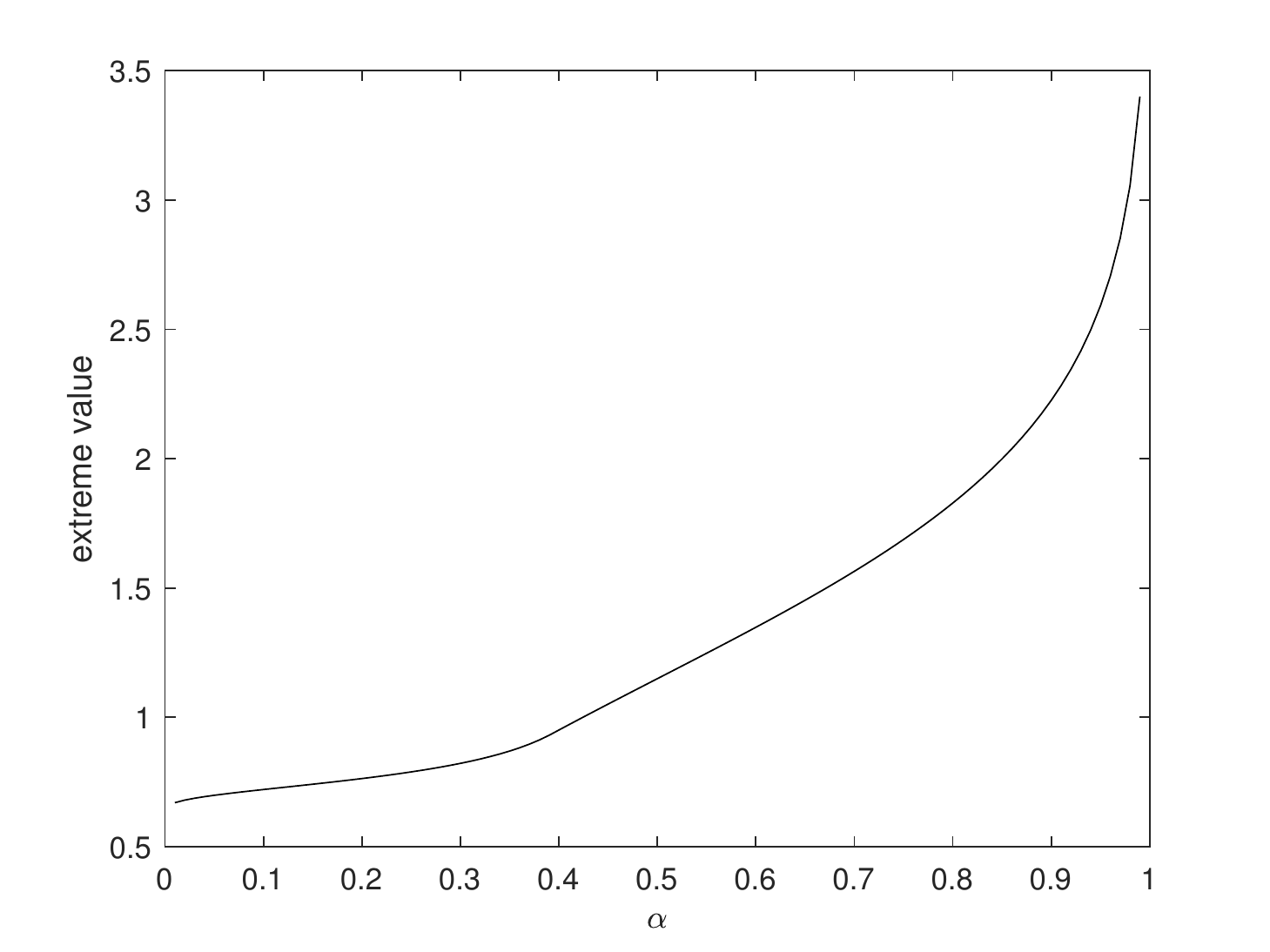}}
	\subfigure[Extreme value error with different order $\nu$]{\label{figure2_b}
		\includegraphics[width=0.48\linewidth]{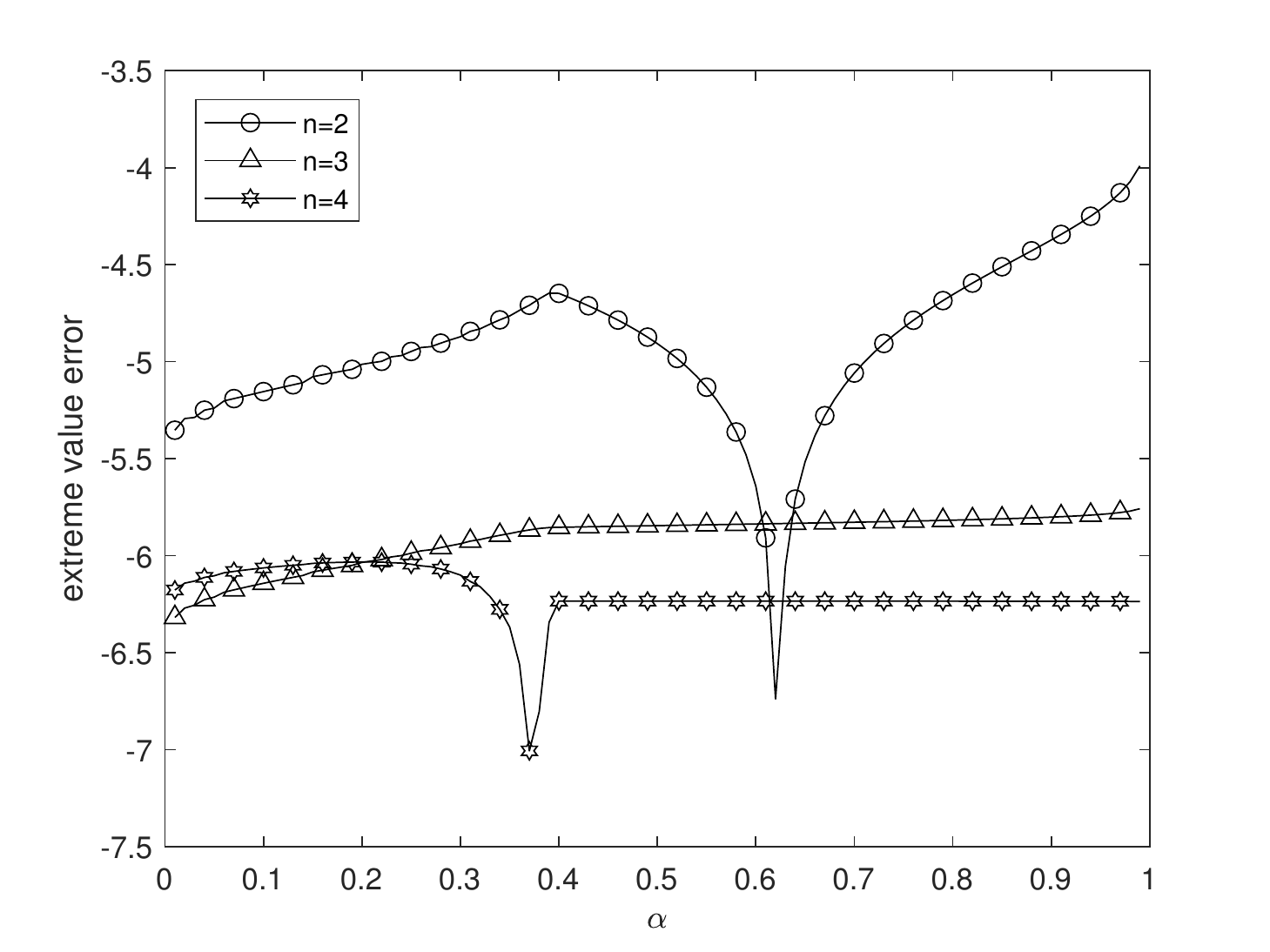}}
	\caption{Plots of the extreme value and error}
\end{figure}

Fig. \ref{figure2_a} shows the extreme value of $\alpha$-path under the given parameters. Fig. \ref{figure2_b} is the image that compares \cref {figure1_a,figure1_b,figure1_c} after logarithmization. As shown in Fig. \ref{figure2_b} and Fig. \ref{figure1sub}, compared with $n\geq 5$, the absolute error of the results is acceptable when $n = 3$ or $4$.

Then we analyze the error in the case of order $\nu\in[0.1,0.9]$ and parameter $\upsilon\in[1,3]$ respectively, and the results are as follows:

\begin{figure}[H]
	\centering  
	\subfigbottomskip=1pt 
	\subfigcapskip=-5pt 
	\subfigure[Error with different order $\nu$]{\label{figure_error1}
		\includegraphics[width=0.48\linewidth]{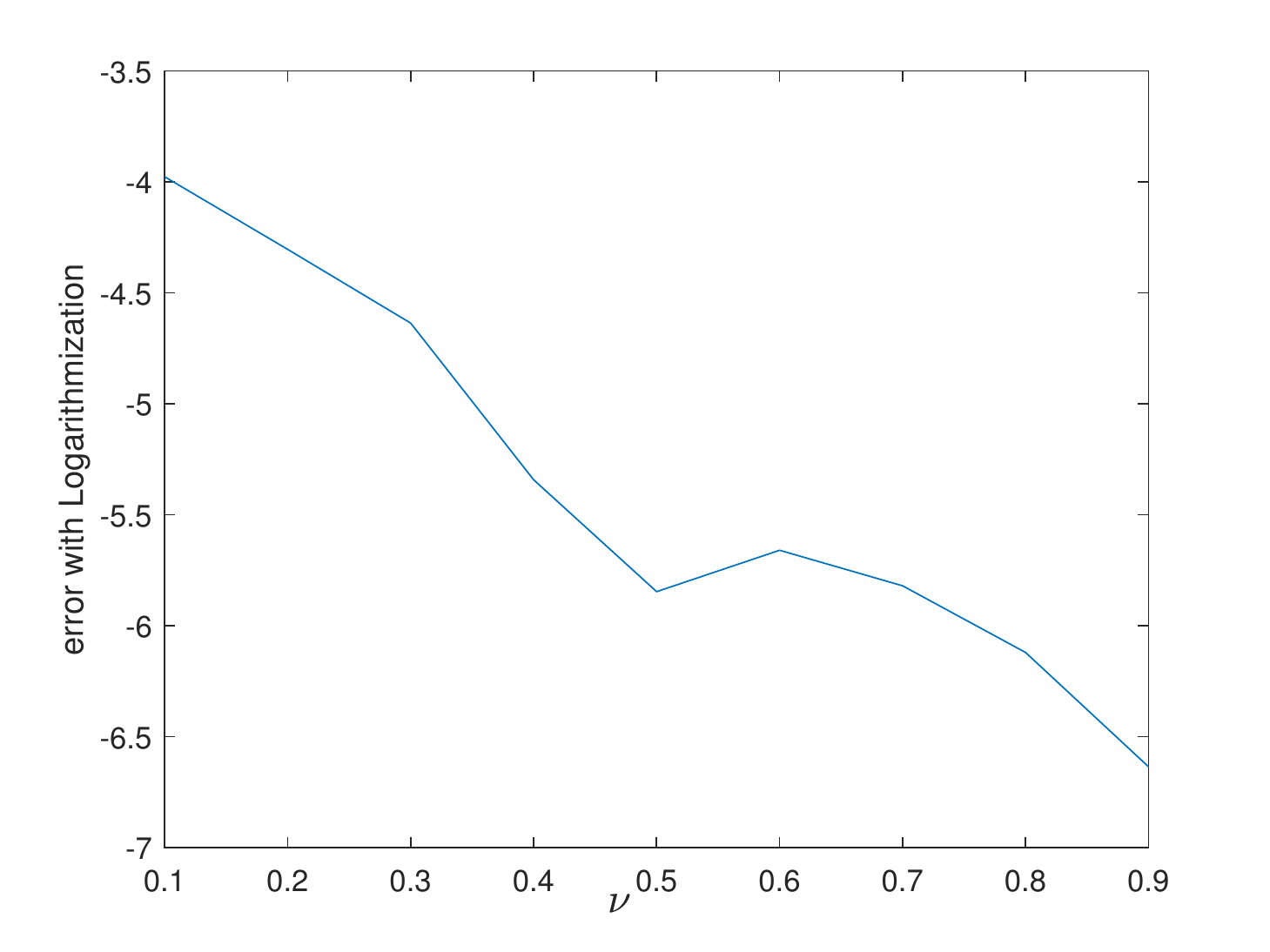}}
	\subfigure[Error with different $\upsilon$]{\label{figure_error2}
		\includegraphics[width=0.48\linewidth]{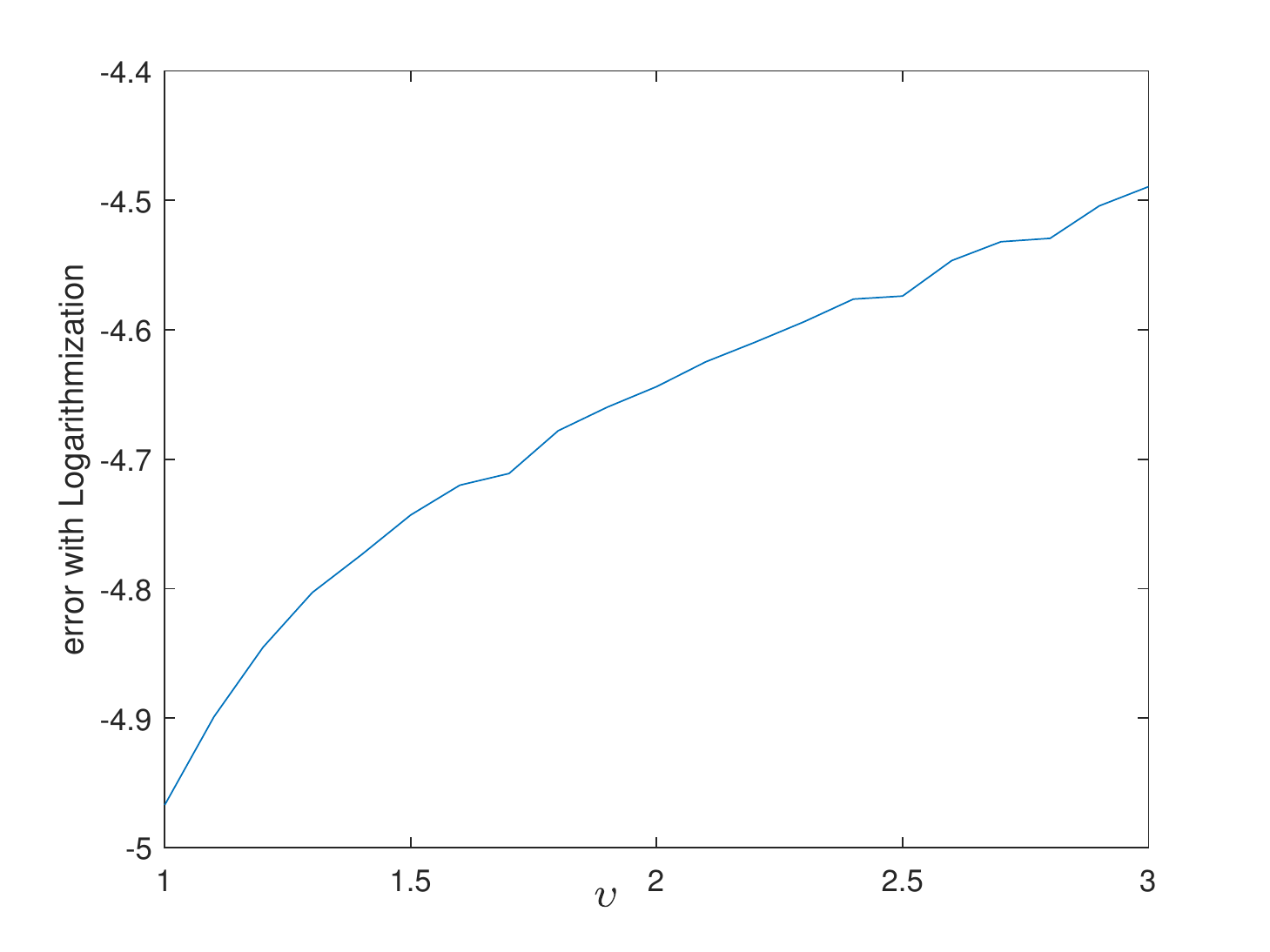}}
	\caption{Plots of the error of different parameters}
\end{figure}

In Fig. \ref{figure_error1}, The $Y$-axis represents the change of MAE after logarithmizing the $alpha$-path, and the $X$-axis represents the change of order of UFDE. In Fig. \ref{figure_error2}, $X$-axis represents the change of parameter $\upsilon$.

Obviously, with the increase of order $\nu$, MAE is also decreasing, which shows that the order of fractional differential equation also has an impact on the error. From Fig. \ref{figure_error1} shows a negative correlation between the two. For the parameter $\upsilon$, it is positively correlated with MAE. With the increase of $\upsilon$, undering the same $t$, the integer derivative value of the right term of the FDE also increases correspondingly. According to the Eq. \eqref{accuracy_pr}, the error bound will become larger.


\subsection{Inverse distribution of UFDE's solution}

In many cases, UFDE is unable to find an exact analytical solution. In this section, we present a numerical method for calculating the inverse uncertain distribution of its solution. The algorithm is given as follows.
\begin{algorithm}[H]
  \caption{Inverse distribution of UFDE's solution.}
  \label{alg:ID_UFDE}
  \begin{algorithmic}[1]
	  \State Write UFDE in the corresponding integral form and split $t$ and $\alpha$ according to the interval $0.01$
	  \State Compute the inverse distribution of UFDE's solution for different $t\in[a, b]$ and $\alpha\in(0,1)$ using Algorithm \ref{alg:Adams}.
  \end{algorithmic}
\end{algorithm}

\begin{figure}[H]

	\centering  
	\subfigbottomskip=1pt 
	\subfigcapskip=-5pt 
	\subfigure[$\alpha$-path($t=1$)]{
		\includegraphics[width=0.48\linewidth]{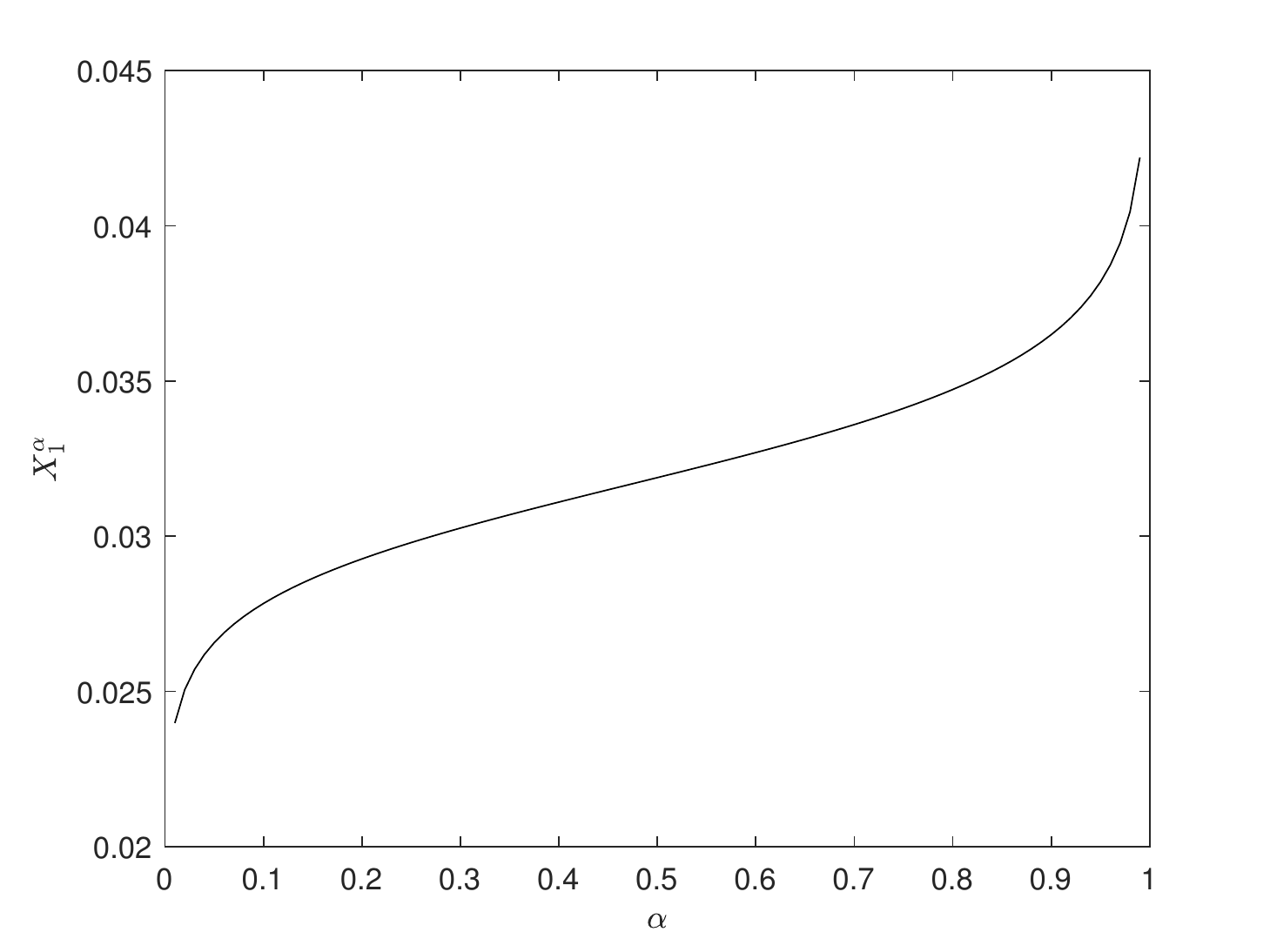}}
	\subfigure[$\alpha=0.1:0.9,t=0.01:1$]{
		\includegraphics[width=0.48\linewidth]{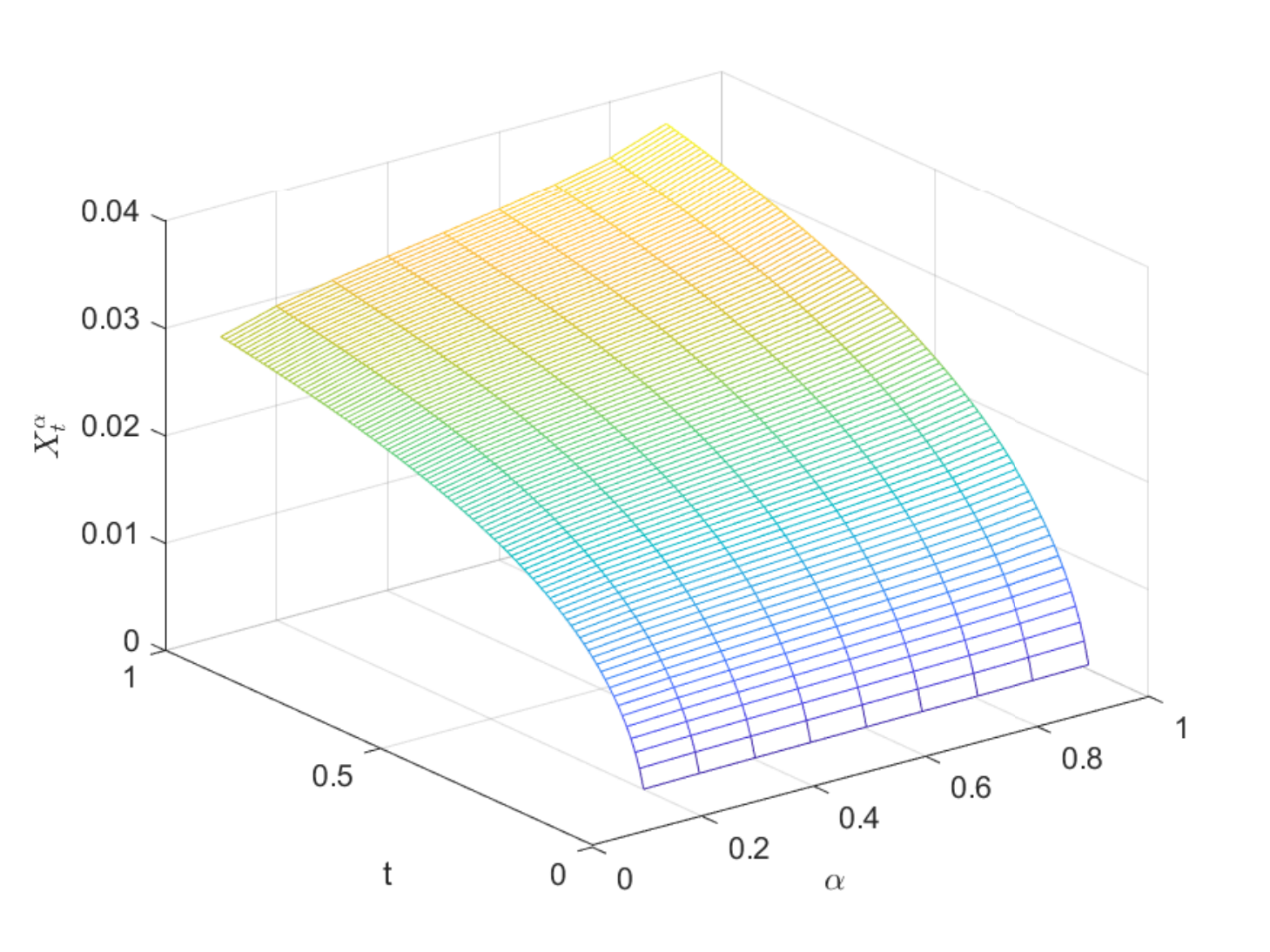}}
	\caption{Numerical results of the uncertain distribution of (\ref{eg2})}
		\label{figure3}
\end{figure}

\begin{example}
Assume the following Caputo type UFDE,
\begin{equation}\label{eg2}
_a^C D_t^{\nu}X_t=a(\mu-X_t)+\sigma\sqrt{X_t}\frac{\mathrm dC_t}{\mathrm dt},
\end{equation}
with initial value $X_0=0$.
\end{example}
When $a=1.2, \mu=0.05, \sigma=0.04, \nu=0.8$, the numerical results of UFDE is shown in Fig. \ref{figure3}.

\subsection{First hitting time of UFDE's solution}
When UFDE is nonlinear, it is generally difficult to find the exact solution. In this subsection, we use example to show how to  apply fractional Adams method to calculate the FHT of nonlinear fractional equation. When $z>J(X_0)$ and $J(x)$ is a nondecreasing function, according to \cite{jin2020first}, the distribution of the FHT is

\begin{equation}\label{eq:fht}
1-\inf\left\{\alpha\in(0,1)\mid\sup\limits_{0\leq t\leq s}J(X_t^{\alpha})\geq z\right\}.
\end{equation}

\begin{algorithm}[H]
  \caption{FHT of UFDE's solution.}
  \label{alg:FHT_UFDE}
  \begin{algorithmic}[1]
	  \State Write UFDE in the corresponding integral form and split $t$ and $\alpha$ according to the interval $0.01$
	  \State Compute the inverse distribution of UFDE's solution for different $t\in[a, b]$ and $\alpha\in(0,1)$ using Algorithm \ref{alg:Adams}.
	  \State Give the FHT for different $\alpha$ according to the Eq. (\ref{eq:fht}).
  \end{algorithmic}
\end{algorithm}

\begin{example}
Assume that a nonlinear UFDE of the Caputo type is
\begin{equation}\label{eg3}
^{C}_aD^{\nu}_tX_t=\sqrt{X_t-1}+(1-t)\frac{\mathrm dC_t}{\mathrm dt},
\end{equation}
with the initial value $X_0=3$ and $J(x)=x$, where $J(X_0)<z$.
\end{example}

\begin{figure}[H]
	\centering  
	\subfigbottomskip=1pt 
	\subfigcapskip=-5pt 
	\subfigure[FHT U(s)]{
		\includegraphics[width=0.48\linewidth]{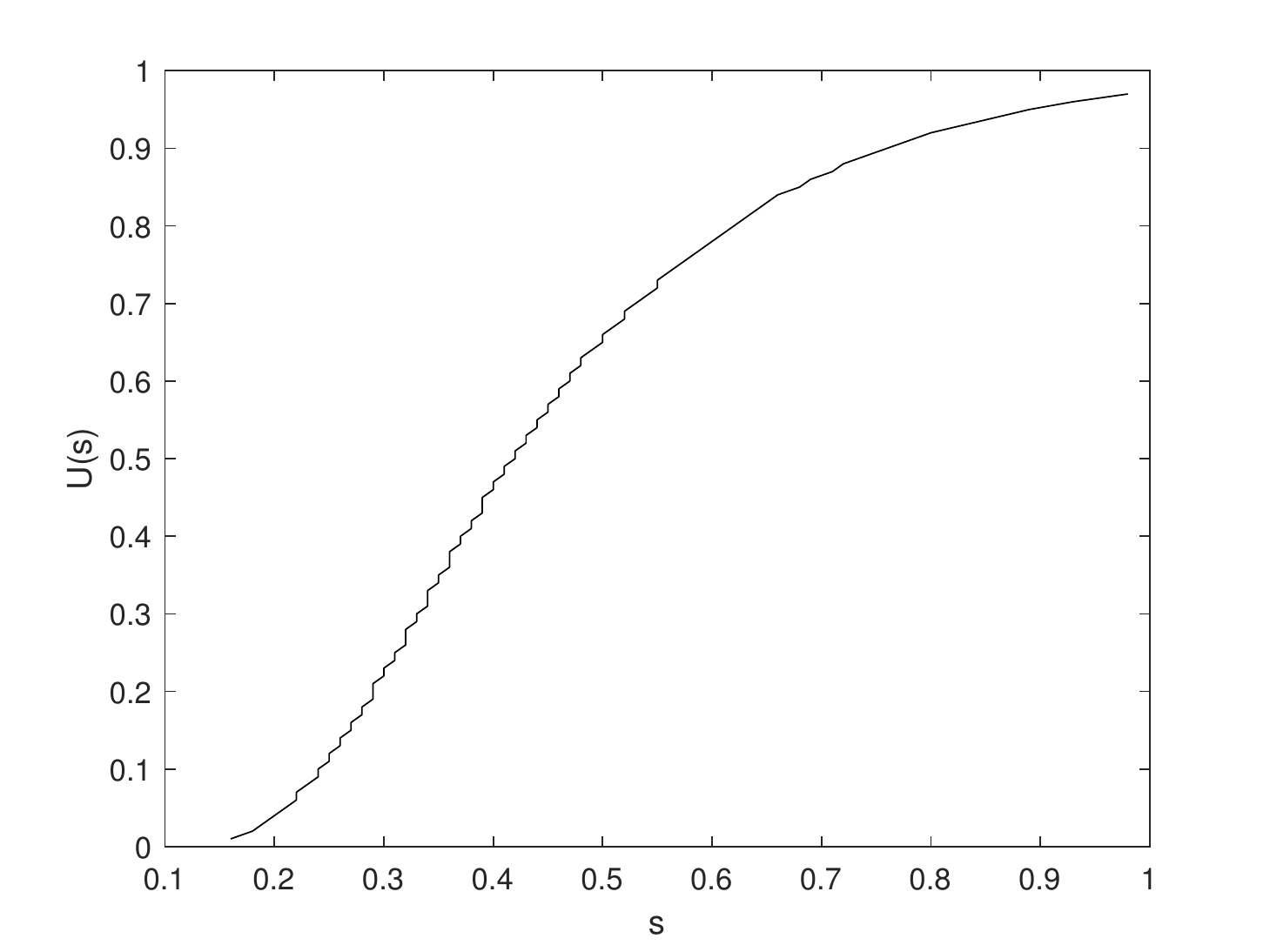}}
	\subfigure[$\alpha$-path($t=1$)]{
		\includegraphics[width=0.48\linewidth]{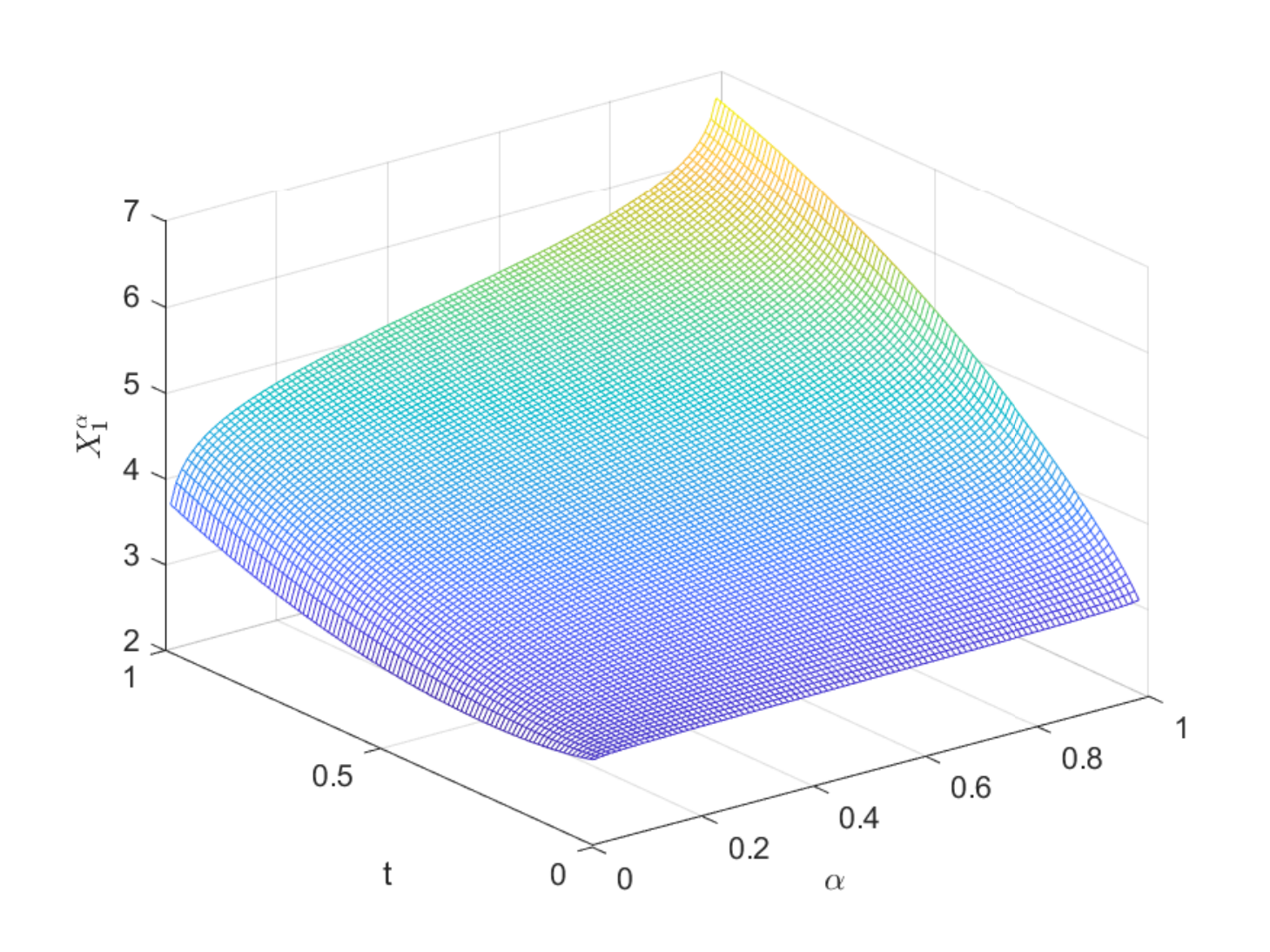}}
	\caption{Results of distribution of the FHT when $J(X_0) < z$.}
	\label{figure4}
\end{figure}

We assign values to the parameters $z=4, \nu=0.8$. According to Algorithm \ref{alg:Adams}, the distribution of FHT is shown in Fig. \ref{figure4}.

In subsection 4.2 and 4.3, it is laborious to verify the accuracy of the results. We only show the rationality of this method through linear and nonlinear examples and calculation of different indicators of the model.
According to the previous equation on accuracy \eqref{accuracy_pr} and the example in subsection 4.1, the accuracy of calculation is guaranteed. Since the calculation of fractional differential equations involves the weighting of the whole time domain, the stability of the algorithm is not given, but only shown from the stability of the numerical results. It can be seen from the numerical simulation that the stability cannot be guaranteed with the increase of $n$.

\section{Conclusion}
A novel numerical method has been developed to solve UFDE with initial value conditions. On the basis of the existing Adams method, according to the product integration method, we gave the polynomial form of lagrange basis function in the case of several values, which extended the traditional Adams method to the fractional Adams method with any node, and did not cause redundancy in computation time. Moreover, three numerical algorithms have been designed for calculating the extreme value, the inverse uncertain distribution and the first hitting time of the solution for UFDE to verify the effectiveness of the methods. Before using the $n$-order Adams method, we need to calculate the first $n$ points, so the accuracy of the initial points is very important. With the aim of more accurate results, we will further study and develop more efficient numerical methods in future research to calculate the estimated values of the first $n$ nodes.

\section*{Data Availability}
Data sharing is not applicable to this article as no datasets were generated or analyzed during the current study.

\section*{Acknowledgments}
This work is supported by the National Natural Science Foundation of China (No.12201304 and No.12071219) and supported by Academic Program Development of Jiangsu Higher Education Institutions (PAPD), Natural Science Foundation of Jiangsu Province (No.BK20210605), the General Research Projects of Philosophy and Social Sciences in Colleges and Universities (2022SJYB0140), the Jiangsu Province Student Innovation Training Program (202110298040Z and 202210298050Z).

\bibliography{bibfile}

\begin{thebibliography}{10}

\bibitem{kahneman2013prospect}
Daniel Kahneman and Amos Tversky.
\newblock Prospect theory: An analysis of decision under risk.
\newblock In {\em Handbook of the fundamentals of financial decision making:
  Part I}, pages 99--127. World Scientific, 2013.

\bibitem{liu2007uncertainty}
Baoding Liu et~al.
\newblock Uncertainty theory.
\newblock In {\em Uncertainty Theory}, pages 205--234. Springer, 2007.

\bibitem{baoding2007uncertainty}
Baoding Liu.
\newblock Uncertainty theory. 2.
\newblock 2007.

\bibitem{liu2012there}
Baoding Liu.
\newblock Why is there a need for uncertainty theory.
\newblock {\em Journal of Uncertain Systems}, 6(1):3--10, 2012.

\bibitem{liu2009some}
Baoding Liu.
\newblock Some research problems in uncertainty theory.
\newblock {\em Journal of Uncertain Systems}, 3(1):3--10, 2009.

\bibitem{liu2008fuzzy}
Baoding Liu.
\newblock Fuzzy process, hybrid process and uncertain process.
\newblock {\em Journal of Uncertain Systems}, 2(1):3--16, 2008.

\bibitem{zhu2010uncertain}
Yuanguo Zhu.
\newblock Uncertain optimal control with application to a portfolio selection
  model.
\newblock {\em Cybernetics and Systems: An International Journal},
  41(7):535--547, 2010.

\bibitem{chen2013uncertain}
Xiaowei Chen, Yuhan Liu, and Ralescu~Dan A.
\newblock Uncertain stock model with periodic dividends.
\newblock {\em Fuzzy Optimization and Decision Making}, 12(1):111--123, 2013.

\bibitem{yang2020parameter}
Xiangfeng Yang, Yuhan Liu, and Park Gyei-Kark.
\newblock Parameter estimation of uncertain differential equation with
  application to financial market.
\newblock {\em Chaos, Solitons \& Fractals}, 139:110026, 2020.

\bibitem{liu2013toward}
Baoding Liu.
\newblock Toward uncertain finance theory.
\newblock {\em Journal of Uncertainty Analysis and Applications}, 1(1):1--15,
  2013.

\bibitem{yao2013numerical}
Kai Yao and Xiaowei Chen.
\newblock A numerical method for solving uncertain differential equations.
\newblock {\em Journal of Intelligent \& Fuzzy Systems}, 25(3):825--832, 2013.

\bibitem{yang2015adams}
Xiangfeng Yang and Ralescu~Dan A.
\newblock Adams method for solving uncertain differential equations.
\newblock {\em Applied Mathematics and Computation}, 270:993--1003, 2015.

\bibitem{gao2016milne}
Rong Gao.
\newblock Milne method for solving uncertain differential equations.
\newblock {\em Applied Mathematics and Computation}, 274:774--785, 2016.

\bibitem{wang2015adams}
Xiao Wang, Yufu Ning, Moughal~Tauqir A, and Xiumei Chen.
\newblock Adams--simpson method for solving uncertain differential equation.
\newblock {\em Applied Mathematics and Computation}, 271:209--219, 2015.

\bibitem{zhang2017hamming}
Yi~Zhang, Jinwu Gao, and Zhiyong Huang.
\newblock Hamming method for solving uncertain differential equations.
\newblock {\em Applied Mathematics and Computation}, 313:331--341, 2017.

\bibitem{Wu2022sest}
Guocheng Wu, Jiali Wei, Cheng Luo, and Lanlan Huang.
\newblock Parameter estimation of fractional uncertain differential equations
  via adams method.
\newblock {\em Nonlinear Analysis: Modelling and Control}, 27:1--15, 02 2022.

\bibitem{Luo_Wu_Huang_2022}
Cheng Luo, Guocheng Wu, and Lanlan Huang.
\newblock Fractional uncertain differential equations with general memory
  effects: Existences and alpha-path solutions.
\newblock {\em Nonlinear Analysis: Modelling and Control}, 28(1):152--179, Dec.
  2022.

\bibitem{lu2019numerical}
Ziqiang Lu and Yuanguo Zhu.
\newblock Numerical approach for solution to an uncertain fractional
  differential equation.
\newblock {\em Applied Mathematics and Computation}, 343:137--148, 2019.

\bibitem{ma2020application}
Yutian Ma and Wenwen Li.
\newblock Application and research of fractional differential equations in
  dynamic analysis of supply chain financial chaotic system.
\newblock {\em Chaos, Solitons \& Fractals}, 130:109417, 2020.

\bibitem{liping2021new}
Liping Chen, Khan~Muhammad Altaf, Atangana Abdon, and Kumar Sunil.
\newblock A new financial chaotic model in atangana-baleanu stochastic
  fractional differential equations.
\newblock {\em Alexandria Engineering Journal}, 60(6):5193--5204, 2021.

\bibitem{kilbas2006theory}
Kilbas~Anatoli{\u\i} Aleksandrovich, Srivastava~Hari M, and Trujillo~Juan J.
\newblock {\em Theory and applications of fractional differential equations},
  volume 204.
\newblock elsevier, 2006.

\bibitem{wu2019new}
Guocheng Wu, ZhenGuo Deng, Baleanu Dumitru, and DeQiang Zeng.
\newblock New variable-order fractional chaotic systems for fast image
  encryption.
\newblock {\em Chaos: An Interdisciplinary Journal of Nonlinear Science},
  29(8):083103, 2019.

\bibitem{abdeljawad2020discrete}
Abdeljawad Thabet, Banerjee Santo, and GuoCheng Wu.
\newblock Discrete tempered fractional calculus for new chaotic systems with
  short memory and image encryption.
\newblock {\em Optik}, 218:163698, 2020.

\bibitem{wu2022unified}
Guocheng Wu, Hua Kong, Maokang Luo, Hui Fu, and Lanlan Huang.
\newblock Unified predictor--corrector method for fractional differential
  equations with general kernel functions.
\newblock {\em Fractional Calculus and Applied Analysis}, 25(2):648--667, 2022.

\bibitem{fu2021fractional}
Hui Fu, Guocheng Wu, Guang Yang, and Lanlan Huang.
\newblock Fractional calculus with exponential memory.
\newblock {\em Chaos: An Interdisciplinary Journal of Nonlinear Science},
  31(3):031103, 2021.

\bibitem{zhu2015uncertain}
Yuanguo Zhu.
\newblock Uncertain fractional differential equations and an interest rate
  model.
\newblock {\em Mathematical Methods in the Applied Sciences},
  38(15):3359--3368, 2015.

\bibitem{zhu2015existence}
Yuanguo Zhu.
\newblock Existence and uniqueness of the solution to uncertain fractional
  differential equation.
\newblock {\em Journal of Uncertainty Analysis and Applications}, 3(1):1--11,
  2015.

\bibitem{jin2019extreme}
Ting Jin, Yun Sun, and Yuanguo Zhu.
\newblock Extreme values for solution to uncertain fractional differential
  equation and application to american option pricing model.
\newblock {\em Physica A: Statistical Mechanics and its Applications},
  534:122357, 2019.

\bibitem{jin2020first}
Ting Jin and Yuanguo Zhu.
\newblock First hitting time about solution for an uncertain fractional
  differential equation and application to an uncertain risk index model.
\newblock {\em Chaos, Solitons \& Fractals}, 137:109836, 2020.

\bibitem{diethelm2002predictor}
Diethelm Kai, Ford~Neville J, and Freed~Alan D.
\newblock A predictor-corrector approach for the numerical solution of
  fractional differential equations.
\newblock {\em Nonlinear Dynamics}, 29(1):3--22, 2002.

\bibitem{li2011numerical}
Changpin Li, An~Chen, and Junjie Ye.
\newblock Numerical approaches to fractional calculus and fractional ordinary
  differential equation.
\newblock {\em Journal of Computational Physics}, 230(9):3352--3368, 2011.

\bibitem{Podlubny1999FractionalDE}
Igor Podlubny.
\newblock Fractional differential equations : an introduction to fractional
  derivatives, fractional differential equations, to methods of their solution
  and some of their applications.
\newblock 1999.

\bibitem{lu2019european}
Ziqiang Lu, Hongyan Yan, and Yuanguo Zhu.
\newblock European option pricing model based on uncertain fractional
  differential equation.
\newblock {\em Fuzzy Optimization and Decision Making}, 18(2):199--217, 2019.

\end{thebibliography}

\end{document}